\documentclass{article}

\usepackage{microtype}
\usepackage{amsfonts}
\usepackage{graphicx}
\usepackage[nice]{nicefrac}
\usepackage{algorithm}% http://ctan.org/pkg/algorithms
\usepackage{algpseudocode}% http://ctan.org/pkg/algorithmicx
\usepackage{hyperref}
\usepackage{amsfonts}
\usepackage{amsthm}
\usepackage{amsmath}
\usepackage{amssymb}
\usepackage{lipsum}
\usepackage{multirow}
\usepackage{authblk}
\usepackage{siunitx}
\usepackage{enumitem}
\usepackage{upgreek}
\usepackage{subfigure}
\usepackage{pgfplots}
\usetikzlibrary{arrows.meta}
\usepackage{xfrac}
\usepackage{etoolbox}
\usepackage{listings}
\usepackage{lmodern}
\usepackage{tabularx}
\usepackage{booktabs}
\pgfplotsset{compat=1.18}
\apptocmd{\thebibliography}{\raggedright}{}{}

\newtheorem{theorem}{Theorem}[section]

\newtheorem{proposition}{Proposition}[section]
\newtheorem{remark}{Remark}[section]
\newtheorem{defn}{Definition}[section]
\newtheorem{eg}{Example}[section]
\numberwithin{equation}{section}

\title{A hierarchical resource-efficient deep policy gradient method for continuous-time optimal control problems}

\author[1]{Arash Fahim\thanks{\href{mailto:arash@math.fsu.edu}{arash@math.fsu.edu}, \url{https://arashfahim.github.io}}}
\author[2]{Md. Arafatur Rahman}

\affil[1]{Department of Mathematics, Florida State University}
\affil[2]{Internal Audit - Model Risk, Citibank N.A., N.Y.}

\usepackage{amsopn}
\usepackage{xcolor}
\selectcolormodel{gray}

\begin{document}

\maketitle

\begin{abstract}
In this paper, we propose an efficient implementation of a deep policy gradient method (PGM) for optimal control problems in continuous time. For continuous-time problems that require a fine time discretization to achieve a desired accuracy, the proposed method improves time efficiency and performance by strategically allocating computational resources across scales, i.e., the number of trajectories, the granularity of time discretization, and the complexity of the neural network architecture. The main idea of the paper is to avoid committing to a fine time discretization. At first, we train a policy, modeled by a neural network, for a discretized optimal control problem in a coarse time scale. Then, we only discretize more in the time intervals where there are indications that the coarse scheme is not accurate enough and train a new policy in the finer scale. We then proceed to refine the time grid further to achieve better accuracy in the new smaller time intervals. Our theoretical result indicates how the new schedule for allocation of resources in different time scales can lead to efficiency. We conclude the paper by numerical experiments on a linear-quadratic stochastic optimal control problem and an optimal execution problem from quantitative finance. 
\end{abstract}

%%%%%%%%%%%%%%%%%%%%%%%%%%%%%%%%
% \begin{list}{-}{Notation}

% \item $\mathcal{N}(\mu,\Sigma)$ is the law of normal distribution with mean $\mu$ and covariance matrix $\Sigma$.
% \end{list}
%%%%%%%%%%%%%%%%%%%%%%%%%%%%%%%%%%%%%
%%%%%%%%%%%%%%%%%%%%%%%%%%%%%%%%%%%%%
\section{Introduction}
\label{sec:intro}

In this paper, we focus on continuous-time stochastic optimal control problems given by
\begin{equation}\label{prob:control}
\begin{split}
    &\inf_{\pi\in\Uppi}\mathfrak{J}(\pi),~\textrm{with}~\mathfrak{J}(\pi):=\mathbb{E}\bigg[\int_0^TL(t,X^\pi_t,\pi_t)dt + g(X_T^{\pi})\bigg]\\
    &dX^\pi_t=\mu(t,X_t^{\pi},\pi_t)dt+\sigma(t,X_t^{\pi},\pi_t)dW_t
\end{split}
\end{equation}
Here $W$ is a Brownian motion, $L$, $g$, $\mu$ and $\sigma$ are Lipschitz continuous functions in all arguments, and $\pi$ is a stochastic process in the set $\Uppi\neq\emptyset$ of admissible policies (controls), which guarantees that \eqref{prob:control} is well defined and $-\infty<\inf_{\pi\in\Uppi}\mathfrak{J}(\pi)$. For a discussion on the admissibility of policies, we refer the reader to \cite[Chapters I, III, and IV]{FS06}. When dealing with a continuous-time problem, we would like to consider a fine time discretization of \eqref{prob:control}:
\begin{align}
        & \inf_{\pi\in\Uppi_\delta}{\mathcal{J}}^n(\pi),~\textrm{with}~{\mathcal{J}}^n(\pi):=\mathbb{E}\bigg[\sum_{i\in[n]}L(t_i,\hat X^{\pi}_{t_i},\pi_{t_i})\delta + g(\hat{X}_T^\pi)\bigg],\label{prob:control_discretized}\\
        &\Delta \hat{X}^\pi_{t_{i+1}}:=\mu(t_i,\hat{X}_{t_i}^\pi,\pi_{t_i}) \delta+ \sigma(t_i,\hat{X}_{t_i}^\pi,\pi_{t_i})\Delta W_{t_{i+1}}\label{Euler-Maruyama}
\end{align}
In the above, $[n]:=\{0,...,n-1\}$, $\delta :=\sfrac{T}{n}$, $t_i=i\delta$, and $\Delta {X}_{t_{i+1}}={X}_{t_{i+1}}-X_{t_i}$ for any continuous-time stochastic process $\{X_t:t\ge0\}$. The admissible class $\Uppi_\delta:=\{\pi\in\Uppi:\pi_s=\pi_{t_i}\text{ a.s.\ for }s\in[t_i,t_{i+1}),\,i\in[n]\}$ consists of the controls in $\Uppi$ that are piecewise constant on the time grid $\{t_i\}$; we identify it, via this piecewise-constant extension, with the set of admissible controls for the discrete-time problem \eqref{prob:control_discretized}, so that a control and its restriction to $\{\pi_{t_i}:i\in[n]\}$ are treated as the same object.
We consider a policy $\pi$ given by a feed-forward neural network; $\pi(\cdot;\theta):[0,T]\times\mathbb{R}^d\to \mathbb{R}^m$, where $\theta$ lies on a high-dimensional Euclidean space and represents the weights and biases of the neural network. Replacing $\pi$ with $\pi(\cdot;\theta)$, we obtain
\begin{equation}\label{prob:discretized-theta}
        \begin{split}
        & \inf_\theta{\mathcal{J}}^n(\theta),~\textrm{with}~{\mathcal{J}}^n(\theta):=\mathbb{E}\bigg[\sum_{i\in[n]}L\big(t_i,\hat X^{\theta}_{t_i}, \pi({t_i},X^{\theta}_{t_i};\theta) \big) \delta + g(\hat{X}_T^\theta)\bigg],\\
        &\Delta \hat{X}^\theta_{t_{i+1}}=\mu\big(t_i,\hat{X}_{t_i}^\theta,\pi(t_i,\hat{X}_{t_i}^\theta;\theta)\big)\delta +\sigma\big(t_i,\hat{X}_{t_i}^\theta, \pi(t_i,\hat{X}_{t_i}^\theta;\theta)\big) \Delta W_{t_{i+1}}
        \end{split}
\end{equation}
Note that the above problem is an optimization problem over variable $\theta$ and not anymore a control problem.
By simulating $M$ independent paths of Brownian motion $W$, $\{W_t^m:t\in[0,T], m\in[M]\}$, we approximate \eqref{prob:discretized-theta} by the following empirical risk minimization problem:
\begin{equation}\label{prob:risk_minimization_}
        \begin{split}
        & \inf_\theta\hat{\mathcal{J}}^n(\theta),~\textrm{with}~\hat{\mathcal{J}}^n(\theta):=\sum_{m\in[M]}\bigg[\sum_{i\in[n]}L\big(t_i,\hat X^{m}_{t_i}, \pi({t_i},X^{m}_{t_i};\theta) \big) \delta + g(\hat{X}_T^m)\bigg],\\
        &\Delta \hat{X}^m_{t_{i+1}}=\mu\big(t_i,\hat{X}_{t_i}^m,\pi(t_i,\hat{X}_{t_i}^m;\theta)\big)\delta +\sigma\big(t_i,\hat{X}_{t_i}^m, \pi(t_i,\hat{X}_{t_i}^m;\theta)\big) \Delta W^m_{t_{i+1}}
        \end{split}
\end{equation}
For $\delta$ significantly small, \eqref{prob:control_discretized} and \eqref{prob:risk_minimization_} are high-frequency discrete-time problems and the implementation of deep numerical methods can lose efficiency in the runtime and memory, especially if the architecture of the neural network is chosen complex to capture variability in time; the number of operations in back-propagation grows linearly in $n$ with a coefficient that depends on the number of sample paths of $\hat{X}$ and the number of parameters of the neural network. In some applications and to a certain extent, the high-frequency nature of the problem dictates accuracy over efficiency. For example, in optimal execution under a stochastic price impact model, one needs to find the solution to a control problem sufficiently quickly to allow for the trading strategy to be adjusted at a high frequency that makes the problem closer to a continuous-time control problem. 

\subsection{Related work}
One of the deep numerical methods for optimal control is the policy gradient method (PGM). PGM approximates an optimal control (policy or action in reinforcement learning) in a class of a semi-parametric functions of time and state variables. PGM evaluates the empirical cost over a collection of non-optimal sample trajectories of the state variable and applies gradient descent method to modify the parameters toward a minimum of the empirical cost function. With modification of parameters, the sample trajectories are also getting closer to the optimal trajectories. PGM for discrete-time optimal control problems was introduced in \cite{BT96} and later in \cite{EH2016}. Application of such methods in quantitative finance was studied in \cite{FMW20, GPW21, GPW22,RS2023, RST23} among other papers. \cite{RS2023} categorized overfitting in terms of convergence properties of empirical value function to the actual value function. Other studies such as \cite{GPW21, GPW22} employed a different neural network for the policy, as a function of state variable only, at each time step and trained the policy stepwise via dynamic programming. This approach does not generalize well in time variable and demands substantial memory allocation when dealing with large number of time steps. In this paper, we follow the approach of \cite{RS2023,RST23} and consider the policy as a single neural network, as a function of time and state variables, to allow generalization in time variable. 

While PGM has been widely studied in discrete time reinforcement learning (RL), e.g., Markov decision processes in \cite{K01,SMSM99}, the application to continuous-time problems has only been considered more recently, e.g. \cite{WZZ20,JZ23,JZ22,GRZ24,ZL23}. \cite{JZ22} developed policy-gradient representations and actor-critic algorithms directly in continuous time and space, providing a general theoretical framework for this line of work. Classic approach to RL classifies RL as an infinite horizon problem. Then, the optimal policy is only a function of the state and the continuous-time and discrete-time versions are only different in the frequency of application of the policy. This differs from a finite horizon control problem with time-variable coefficients. In the recent formulation of the RL problem in continuous time with finite horizon, \cite{WZZ20} considered an optimal control problem with the addition of an entropy term to the loss function. This approach guarantees the randomization of the policy, which is a crucial component in RL. One notable implication of this formulation is that, with the addition of an entropy term, a nonconvex control problem becomes an RL problem convex in the randomized policy. This guarantees the convergence of numerical methods to the global minimum. In \cite{GRZ24}, the authors cited \cite{M06, PKK21} as empirical studies that show the discrete-time approximation of continuous-time RL problems of \cite{WZZ20} exhibits performance degradation due to difference in the scale. They showed, theoretically and numerically, that if the gradient is rescaled by $\delta t$, the length of discrete time mesh, then gradient descent method converges in smaller number of epochs. That is a different issue from what we are tackling in this paper.

For optimal control, the randomization is not present. Therefore, we lose the convexity, and simple convergence results are not available except if the control problem itself is convex. See \cite{ZL23} for convergence of PGM for optimal control problems under assumptions of smoothness and convexity. \cite{RSZ23} established linear convergence of a policy-gradient method to stationary points for a class of finite-horizon stochastic control problems, including nonlinear diffusions with degenerate noise. \cite{HXY21} proved convergence of policy-gradient methods for the finite-horizon noisy linear-quadratic regulator, with an optimal-liquidation example closely related to the optimal execution problem we consider in Section~\ref{sec:numerics}. For non-convex Hamiltonian, the question of convergence is analogous to the convergence of gradient descent in non-convex optimization and should be addressed separately. For example, see \cite{RRT17,ZLC17,CCDPW20,DFGLL20} for treatment of non-convexity in gradient methods driven by a Langevin dynamics.

For the continuous-time optimal control problems, an approximation method is coming from the numerical schemes for the corresponding parabolic or elliptic Hamilton-Jacobi-Bellman equations (HJB). In low dimensions, classical methods such as finite difference are effective in solving HJB equations. For example, see \cite{LM80,TZ97,BZ03}. In higher dimensions, studies \cite{Z01,BT04,LGW06,FTW11,BF14,ZZ14} addressed Monte Carlo numerical schemes to nonlinear PDEs arising in stochastic optimal control via backward stochastic differential equation (BSDE). 
More recently, \cite{JHW17,JHW18,BEJ19} used the BSDE representation of the HJB as a basis for numerical approximation of the value function via deep neural networks. Another approach toward solving the HJB equation is the deep Galerkin method in \cite{SS18}, which is related to the Physics informed neural networks.

For deterministic optimal control problems, the Bellman equation is not parabolic anymore and the solution may not be regular. This may cause issues in stability of numerical schemes due to irregularity of the solution. In high dimensions, the only viable option to numerically solve the Bellman equation is the Hopf-Lax method. This method relies on some restrictive assumption on the problem and uses ODEs from the characteristics to overcome the curse of dimensionality. More specifically, it relies on generating trajectories of the primal and dual variables forward in time and optimizing with respect to the starting position of the dual trajectory. Generally, curse of dimensionality is overcome by methods that rely on the trajectories rather than discretization of space. For more information on the Hopf-Lax formula in deterministic control, see \cite{DO16,CDOS17,YD21}. While the Hopf-Lax method is efficient in evaluation of a single optimal trajectory, it has not yet been used for training a deterministic optimal control policy for a class of trajectories at different starting positions or starting times.

\subsection{Our contribution}
In this paper, we propose a hierarchical implementation of the PGM for continuous-time optimal control problems, which allows for a systematic management of allocation of computational resources. The computational resources considered in this study include the number of parameters in the neural network and the total number of sample trajectories used in training. Our method is particularly effective when the dynamics are non-stationary in time, i.e., when the state's behavior can vary significantly within a single coarse time step.
Our proposed method, \emph{deep hierarchical PGM}, consists of multiple steps. At each step of the algorithm, the time discretization becomes finer based on the information obtained at the previous step. In the first step, we choose a coarse time-step for problem \eqref{prob:control_discretized}, \emph{coarse problem}, and use PGM on \eqref{prob:risk_minimization_} to train a neural network for the \emph{optimal coarse policy} as a function of time and state variables. After training the optimal coarse policy, we generate a sample of optimal coarse trajectories to find data points at the beginning and end of each coarse time interval and evaluate the cost along each trajectory. Furthermore, we evaluate the \emph{coarse value function} based on the data points on the trajectories and their associated cost. We then compare the distributions of data points and the coarse value functions at the beginning and the end of each coarse interval. The intervals that show the largest difference between the start time and end time are chosen for the next scale.
The points on the coarse trajectory at the beginning of the chosen intervals serve as the initial positions to generate trajectories of the state variable inside the interval on a finer time scale. For each selected interval, we discretize the control problem to a finer time discretization and introduce a new neural network to train an \emph{optimal fine policy}. In the refined control problem, the running cost is the same, but the terminal cost is evaluated by the value function at the end of the interval. We combine the discretized loss functions for all the selected intervals into one single loss function to train the optimal fine policy. We can continue the hierarchical PGM method to finer intervals by choosing a new set of intervals in the fine scale and by training a new value function on the finer intervals.
Breaking the learning problem into multiple stages of refining time discretization allows for flexibility on the number of samples or the architecture of the policy at each scale. We leverage this property to improve the efficiency of the PGM. Our main theoretical result in this paper, Theorem~\ref{thm:management}, shows how different choices of resource allocation contribute to efficiency. The benchmark to compare efficiency with is the \emph{brute-force} implementation of PGM which uses the finest time step in our deep hierarchical PGM across the whole horizon of the problem. In our numerical experiments, we applied the method to two problems, a standard stochastic linear-quadratic optimal control problem and an optimal execution problem under stochastic price impact and stochastic resilience. The first numerical experiment, Section~\ref{sec:lqsc}, is on a 1-dimensional stochastic control problem and the evaluation of the value function is straightforward. The second numerical experiment, Section~\ref{sec:oe}, is on a 4-dimensional optimal execution problem with a mixture of stochastic and deterministic state variables. The degeneracy caused by the deterministic variables is overcome by adding a retraining sub-step in the training of the surrogate for the coarse value function. Additionally, the total-liquidation constraint in the optimal execution problem requires us to include additional trajectories in training of the finer policy. In both cases, the hierarchical PGM demonstrates a gain in efficiency as shown in Section~\ref{sec:numerics}, namely Figures~\ref{fig:value0_comparison_bw_brute_multi} and \ref{fig:value2_comparison_bw_brute_multi} and Tables~\ref{tab:training-times}, \ref{tab:wallclock_lqc_v9}, \ref{tab:cost_balance_stat}, and \ref{tab:wallclock_lob_v7}.

This paper is organized as follows. To make the paper accessible to a broader audience, we provide a brief review of PGM and a discussion on the computational cost of the deep PGM in Section~\ref{sec:prelim}. Section~\ref{sec:main} covers the main contribution of this paper by presenting the deep hierarchical PGM and the main result of the paper, Theorem~\ref{thm:management}. The last section lays out the implementation details for our numerical experiments. For the sake of completeness, the appendix provides a strong error estimate for the discretization of continuous-time stochastic control problems.

%%%%%%%%%%%%%%%%%%%%%%%%%%%%%%%%%%%%%
%%%%%%%%%%%%%%%%%%%%%%%%%%%%%%%%%%%%%
\section{Preliminaries}
\label{sec:prelim}
Throughout this paper, we use the notation $[n]:= \{0,...,n-1\}$, which leads to $[n+1]=\{0,...,n\}$, $[n]+1=\{1,...,n\}$, and $a[n]=\{0,a,...,a(n-1)\}$.

\subsection{PGM for discrete-time optimal control problems}
 \label{sec:pgm_discrete}
For pedagogical reasons, we first review PGM in a general discrete-time setting; readers already familiar with PGM may skip directly to Section~\ref{sec:pgm_conti}. Consider the discrete-time control problems
\begin{equation}\label{prob:discrete_soc}
\begin{split}
  &\inf_{\pi\in\Uppi} \mathcal{J}(\pi),\text{ with } \mathcal{J}(\pi):=\mathbb{E}\bigg[\sum_{i\in[n]}L({i},X_{i},\pi_{i})+g(X_n)\bigg]\\
& X_{i+1} = X_{i}+f({i},X_{i},\pi_{i},\omega_{i+1})
\end{split}
\end{equation}
Here $\Uppi$ is the set of all \emph{admissible} controls, $L:[n]\times\mathbb{R}^d\times\mathbb{R}^m\to\mathbb{R}$ is the running cost, $g:\mathbb{R}^d\to\mathbb{R}$ is the terminal cost, $f:[n]\times\mathbb{R}^d\times\mathbb{R}^m\times\Omega\to\mathbb{R}^d$ is the dynamics of $X$, $(\Omega,\mathbb{P})$ is a probability space, $\mathbb{E}$ is the expectation, and $\{\omega_{i+1}\}_{i\in[n]}$ is a sequence of i.i.d. random variables. The \emph{policy gradient method} (PGM) models the control $\pi_{i}$ in \eqref{prob:discrete_soc} by a parametrized feedback control, $\phi_i(x;\uptheta)$ and reduces it to the following risk minimization problem:
\begin{equation}\label{prob:risk_minimization_theta_discrete}
\begin{split}
  &\inf_{\theta} \mathbb{E}\bigg[\sum_{i\in[n]}L({i},X_{i},\phi_{i}(X_{i};\theta))+g(X_n)\bigg]\\
& X_{{i+1}} = X_{i}+f({i},X_{i},\phi_{i}(X_{i};\theta),\omega_{i+1})
\end{split}
\end{equation}
By sampling $\{\omega_{i+1}\}_{i\in[n]}$ and simulating paths on $X$, we obtain the empirical risk minimization problem:
\begin{equation}\label{prob:empirical_risk_minimization_theta_discrete}
  \begin{split}
  &\inf_{\theta} \sum_{m\in[M]}\sum_{i\in[n]}L\big({i},X^m_{i},\phi_{i}(X^m_{i};\theta)\big)+g(X^m_n)\\
& X^m_{i+1} = X^m_{i}+f\big({i},X^m_{i},\phi_{i}(X^m_{i};\theta),\omega^m_{i+1}\big)
\end{split}
\end{equation}
where $\{\omega^m_{i+1}:i\in[n],m\in[M]\}$ are i.i.d. samples of $\{\omega_{i+1}\}_{i\in[n]}$.
\subsection{PGM for continuous-time optimal control problems}
 \label{sec:pgm_conti}
To apply policy gradient method to the continuous-time problem \eqref{prob:control}, we first discretize time by $\delta:= \sfrac{T}{n}$, $t_i=i\delta$, and $\Delta W_{t_{i+1}}:=W_{t_{i+1}}-W_{t_{i}}$ to obtain the discrete-time problem \eqref{prob:control_discretized}.
% \begin{align}
%     \inf_{\pi\in\Uppi}\mathcal{J}^n(\pi) &,~\text{with}~\mathcal{J}^n(\pi):=\mathbb{E}\left[\sum_{i\in[n]}L(t_{i},\hat{X}^\pi_{t_{i}},\pi_{t_{i}})\delta + g(\hat{X}_T^\pi)\right] \label{prob:control_discretized}\\
%    & \hat{X}^\pi_{t_{i+1}}=\hat{X}^\pi_{t_{i}}+\mu({t_{i}},\hat{X}_{t_{i}}^\pi,\pi_{t_i})\delta+\sigma({t_{i}},\hat{X}_{t_{i}}^\pi,\pi_{t_{i}})\delta W_{{t_{i+1}}}\label{Euler-Maruyama}
% \end{align}
% \begin{equation}\label{prob:control_discretized}
%         \begin{split}
%         & \inf_\pi{\mathcal{J}}^n(\pi),~\textrm{with}~{\mathcal{J}}^n(\pi):=\mathbb{E}\bigg[\sum_{i\in[n]}L(t_i,\hat X^{\pi}_{t_i},\pi_{t_i})\delta + g(\hat{X}_T^\pi)\bigg],\\
%         &\delta \hat{X}^\pi_{t_{i+1}}:=\mu(t_i,\hat{X}_{t_i}^\pi,\pi_{t_i}) \delta+ \sigma(t_i,\hat{X}_{t_i}^\pi,\pi_{t_i})\delta W_{t_{i+1}}
%         \end{split}
% \end{equation}
Given that an optimal control for \eqref{prob:control_discretized} exists, $\pi^{*n} \in\mathop{\text{\normalfont argmin}}_{\pi\in\Uppi_\delta}\mathcal{J}^n(\pi)$,
the discrete-time value function for \eqref{prob:control_discretized} at $\{t_{i}:i\in[n+1]\}$ is given by
\begin{equation}\label{defn:value_discrete}
    \hat{V}(t_{i},x):=\mathbb{E}\Bigg[\sum_{\hat{i}\in[n-i]}L(t_{i+\hat{i}},\hat{X}^{\pi^{*n}}_{t_{i+\hat{i}}},\pi^{*n}_{t_{i+\hat{i}}})\delta + g(\hat{X}_T^{\pi^{*n}})\bigg|\hat{X}^{*n}_{t_i}=x\Bigg], ~~\hat{V}(T,x)=g(x)
\end{equation}
The discretization error can be estimated via the following theorem. The proof of this classical result is provided in Appendix~\ref{sec:error} for completion.
\begin{theorem}\label{thm:strong_error}
    Under the assumptions \textbf{A.1}--\textbf{A.2} in Section~\ref{sec:error}, there exists a $C>0$ independent of $N$ such that
    \begin{equation}
        -C\delta^{\sfrac{1}{4}}\le\sup_{x,i\in[N]}\big(V(t_{i},x)-\hat{V}(t_{i},x)\big)\le C\delta^{\sfrac{1}{2}}
    \end{equation}
    where $V(t,x)$ is the value function of the continuous-time control problem \eqref{prob:control} given by
    \begin{equation}\label{value_fnc_continuous}
        V(t,x):= \inf_{\pi\in\Uppi_t}\mathbb{E}\left[\int_t^TL(s,X^\pi_s,\pi_s)ds + g(X_T^\pi)\bigg| X^\pi_t=x\right]
    \end{equation}
    where $\Uppi_t$ is the set of admissible controls restricted on $[t,T]$. Constant $C$ depends only on $T$ and Lipschitz constant for $\mu$, $\sigma$, $L$, and $g$. 
\end{theorem}
To apply PGM to \eqref{prob:control}, we shall apply it to the discretized version \eqref{prob:control_discretized} through solving the risk minimization problem \eqref{prob:empirical_risk_minimization_theta}. To do so, we need to simulate sample paths of the Euler-Maruyama discretization \eqref{Euler-Maruyama}, $\{\hat{X}_{t_{i-1}}^m:i\in[n+1],m\in[M]\}$, based on samples of Brownian motion $\{W^m_{t_{i}}:i\in[n],m\in[M]\}$ and the samples of $\hat{X}_0$ are drawn from a given initial distribution $\mathcal{D}_0$. 

%%%%%%%%%%%%%%%%%%%%%%%%%%%%%%%%%%%%%
\subsection{Dynamic programming principle and generalization of policy}
Within the method, we locally apply dynamic programming principle (DPP) to generalize the policy inside a time interval. Recall the definition of value function $V(t,x)$ in \eqref{value_fnc_continuous} for the continuous-time control problem \eqref{prob:control}. The simplest form of DPP asserts that for $0\le t<s\le T$, we have
  \begin{equation}\label{DPP_continuous}
        V(t,x)= \inf_{\pi\in\Uppi_{t,s}}\mathbb{E}\left[\int_t^sL(r,X^\pi_r,\pi_r)dr + V(s,X_s^\pi)\bigg| X_t=x\right]
    \end{equation}
    where $\Uppi_{t,s}$ is the set of admissible controls restricted on $[t,s]$. 
    
    Given $V(s,x)$ is known or at least approximated, we apply PGM to minimize right-hand side of \eqref{DPP_continuous} to generalize the approximation of an optimal policy over the interval $[t,s]$. More precisely and similar to (\ref{prob:control_discretized}-\ref{Euler-Maruyama}), we set $\delta=\sfrac{(s-t)}{n}$, $r_i=t+i\delta$, and $\Delta W_{{r_{i}}}=W_{{r_{i}}}-W_{{r_{i-1}}}$, and solve the following discrete-time control problem
 \begin{equation}
 \begin{split}\label{prob:risk_minimization_theta_t_s}
    \inf_{\pi\in\Uppi_{t,s,\delta}}~ &\mathbb{E}\left[\sum_{i\in[n]}L(r_{i},\hat{X}^\pi_{r_{i}},\pi(r_{i},\hat{X}^\pi_{r_{i}}))\delta + V(s,\hat{X}_s^\pi)\Big| \hat{X}^\pi_{t}=x \right] \\
   &\hat{X}^\pi_{r_{i+1}} =\hat{X}^\pi_{r_{i}}+\mu({r_{i}},\hat{X}_{r_{i}}^\pi,\pi(r_{i},\hat{X}^\pi_{r_{i}}))\delta+\sigma({r_{i}},\hat{X}_{r_{i}}^\pi,\pi(r_{i},\hat{X}^\pi_{r_{i}}))\Delta W_{{r_{i+1}}}
\end{split}
\end{equation}
As in \eqref{prob:control_discretized}, $\Uppi_{t,s,\delta}:=\{\pi\in\Uppi_{t,s}:\pi_r=\pi_{r_i}\text{ a.s.\ for }r\in[r_i,r_{i+1}),\,i\in[n]\}$ denotes the controls in $\Uppi_{t,s}$ piecewise constant on the grid $\{r_i\}$, identified with the set of admissible controls for the discrete-time problem \eqref{prob:risk_minimization_theta_t_s}.

The following classical result shows that finding an optimal control for a sufficiently fine discrete-time problem \eqref{prob:risk_minimization_theta_t_s}, with approximated value function $\hat{V}$ in place of $V$, yields an approximately optimal control for the continuous-time problem \eqref{DPP_continuous}. A proof of this result is provided in Section~\ref{sec:error}.
\begin{proposition}\label{prop:epsilon_optimal}
Fix $0\le t_n <t_m \le T$ and let Assumptions \textbf{A.1}--\textbf{A.2} in Section~\ref{sec:error} hold and $\hat{V}$ be an approximation of the value function $V$ such that 
\begin{equation}
    \sup_x|V(t_m,x)-\hat{V}(t_m,x)|\le \dfrac{\epsilon}{2}
\end{equation}
 Define $\pi^\epsilon_s:=\pi^{*n}(t_j,\hat{X}^{\pi^{*n}}_{t_j})$ on $s\in[t_j,t_{j+1})$ where $\pi^{*n}$ is an admissible discrete-time Markovian optimal control such that
\begin{equation}
    \pi^{*n} \in\mathop{\text{\normalfont argmin}}_{\pi\in\Uppi_{t_n,t_m,\delta}}\mathbb{E}\left[\sum_{i=n}^{m-1}L(t_{i},\hat{X}^\pi_{t_{i}},\pi(t_{i},\hat{X}^\pi_{t_{i}}))\delta + \hat{V}(t_m,\hat{X}_{t_m}^\pi)\Big| \hat{X}^\pi_{t}=x\right]
\end{equation}
Then, for sufficiently small $\delta>0$, $\pi^\epsilon_s$ is $\epsilon$-optimal for
\begin{equation}
\inf_{\pi\in\Uppi_{t_n,t_m}}\mathbb{E}\left[\int_{t_n}^{t_m}L(r,X^\pi_r,\pi_r)dr + V(t_m,X_{t_m}^\pi)\bigg| X^\pi_t=x\right]
    \end{equation}
\end{proposition}

%%%%%%%%%%%%%%%%%%%%%%%%%%%%%%%%%%%%%
\subsection{Deep PGM}
Similar to \eqref{prob:risk_minimization_theta_discrete}, in \eqref{prob:control_discretized}, we replace $\Uppi_\delta$ by a parametrized class of controls. More precisely, $\Uppi_\delta=\{\phi(t,x;\theta):\theta\in \mathbb{R}^q\}$, where
\begin{equation}\label{eqn:deepnn}
    \phi(t,x;\theta)=\Big(L_\ell\circ \Upsigma \circ\cdots \circ \Upsigma \circ L_1\Big)(t,x)
\end{equation}
where $L_k(v) = \mathbf{W}_kv+\mathbf{b}_k$ is an affine function from $\mathbb{R}^{k_i}\to \mathbb{R}^{k_{i+1}}$ with $k_1=d+1$, $k_\ell=m$, $\Upsigma$ is an activation function such as ReLU, sigmoid, $\tanh$, and the like, and $\theta=(\mathbf{W}_1,\mathbf{b}_1,...,\mathbf{W}_\ell,\mathbf{b}_\ell)\in\mathbb{R}^q$ with $q = (d+2) k_2 + (k_2+1)k_3+\cdots (k_{\ell-1}+1)k_\ell$.
In problem \eqref{prob:control_discretized}, $\Uppi_\delta$ is replaced by $\{\phi(t,x;\theta):\theta\in\mathbb{R}^q\}$ to yield the following optimization problem:
\begin{align}
    &\inf_{\theta}\mathcal{J}^{n}(\theta),~\mathcal{J}^{n}(\theta):=\mathbb{E}\left[\sum_{i\in[n]}L(t_{i},\hat{X}^\theta_{t_{i}},\phi(t_{i},\hat{X}^\theta_{t_{i}};\theta))\delta + g(\hat{X}_T^\theta)\right] \label{prob:risk_minimization_theta}\\
   & \hat{X}^\theta_{t_{i+1}}=\hat{X}^\theta_{t_{i}}+\mu({t_{i}},\hat{X}_{t_{i}}^\theta,\phi(t_{i},\hat{X}^\theta_{t_{i}};\theta))\delta\nonumber\\
   &\hspace{3cm}+\sigma({t_{i}},\hat{X}_{t_{i}}^\theta,\phi(t_{i},\hat{X}^\theta_{t_{i}};\theta))\Delta W_{{t_{i+1}}},~\hat{X}_0^\theta\sim\mathcal{D}_0\label{X^theta}
\end{align}
In the above, $\mathcal{D}_0$ is an arbitrary distribution with support on a region of interest.
The above minimization can be approximated by an empirical risk minimization problem via simulating the sample paths $\{\hat{X}_{t_{i-1}}^{m,\theta}:m\in[M], i\in[n+1]\}$ of 
 \eqref{X^theta}:
\begin{equation}\label{prob:empirical_risk_minimization_theta}
    \inf_{\theta}\hat{\mathcal{J}}^{n}(\theta),~\hat{\mathcal{J}}^{n}(\theta):=\sum_{m\in[M]}\left[\sum_{i\in[n]}L(t_{i},\hat{X}^{m,\theta}_{t_{i}},\phi(t_{i},\hat{X}^{m,\theta}_{t_{i}};\theta))\delta + g(\hat{X}_T^{m,\theta})\right]
\end{equation}

\begin{defn}[Empirical training domain]\label{defn:empirical_training_domain}
For $i=0,\ldots,N$, the \emph{empirical training domain at time $t_i$} of the value function is the collection of sample points
\begin{equation}
    D_{\text{\normalfont tr}}(t_i):=\{\hat{X}^{m,\theta^{*}}_{t_{i}}:m\in[M]\}
\end{equation}
used in \eqref{value_fnc:risk_minimization} to train the surrogate value function $\chi(t_i,\cdot;\rho)$. The \emph{empirical training domain} of the value function is $D_{\text{\normalfont tr}}:=\bigcup_{i=0}^{N}D_{\text{\normalfont tr}}(t_i)$.
\end{defn}
%%%%%%%%%%%%%%%%%%%%%%%%%%%%%%%%%%%%%%%%%%%%%%%%
%%%%%%%%%%%%%%%%%%%%%%%%%%%%%%%%%%%%%%%%%%%%%%%%
%%%%%%%%%%%%%%%%%%%%%%%%%%%%%%%%%%%%%%%%%%%%%%%%
\subsection{Surrogate for value function}
The goal of PGM is to find an approximately optimal policy rather than the value function. However, to develop the methods of this paper, we need to find an approximate value function for the discretized problem, which is covered in this section. Assume that the optimal policy for \eqref{prob:control_discretized} exists and let $\hat{X}^{\pi^{*n}}$ be an optimal trajectory for discretized-state variable in \eqref{Euler-Maruyama}, we recall from \eqref{defn:value_discrete} that 
$\hat{V}(t_i,x)=\mathbb{E}[Y^{\pi^{*n}}_{t_i}|\hat{X}^{\pi^{*n}}_{t_i}=x]$, where  
\begin{equation}
    Y^{\pi^{*n}}_{t_i}:=\sum_{\hat{i}\in[n-i]}L(t_{i+\hat{i}},\hat{X}^{\pi^{*n}}_{t_{i+\hat{i}}},\pi^{*n}_{t_{i+\hat{i}}})\delta + g(\hat{X}_T^{\pi^{*n}})
\end{equation}
is the cost along the optimal trajectory $\hat{X}^{\pi^{*n}}$. 
On the other hand, by definition of conditional expectation, the value function $\hat{V}(t_i,x)$ is the a.s. unique minimizer of 
\begin{equation}\label{cond_expect}
    \mathop{\inf}\limits_{v(x)}\mathbb{E}[(v(\hat{X}^{\pi^{*n}}_{t_i})-Y^{\pi^{*n}}_{t_i})^2]
\end{equation}
where the infimum is over the set of all Borel functions $v(x)$ such that $v(\hat{X}^{\pi^{*n}}_{t_i})\in L^2$, i.e., $\mathbb{E}[(v(\hat{X}^{\pi^{*n}}_{t_i}))^2]<\infty$.
When the discretized problem is replaced by the risk minimization \eqref{prob:empirical_risk_minimization_theta}, we use an approximate optimal trajectory to approximate the value function.
Given $\theta^*\in\mathop{\text{\normalfont argmin}}\limits_{\theta}{\mathcal{J}}^{n}(\theta)$, \eqref{cond_expect} is replaced by 
\begin{equation}\label{cond_expect_theta}
    \mathop{\min}\limits_{v(x)}\mathbb{E}[(v(\hat{X}^{\theta^{*}}_{t_i})-Y^{\theta^*n}_{t_i})^2]
\end{equation}
where $\hat{X}^{\theta^*}$ satisfies \eqref{X^theta} and
\begin{equation}
    Y^{\theta^{*}}_{t_i}:=\sum_{\hat{i}\in[n-i]}L\big(t_{i+\hat{i}},\hat{X}^{\theta^{*}}_{t_{i+\hat{i}}},\phi(t_{i+\hat{i}},\hat{X}^{\theta^*}_{t_{i+\hat{i}}};\theta^*)\big)\delta + g(\hat{X}_T^{\theta^{*}})
\end{equation}
is the cost along the optimal trajectory $\hat{X}^{\theta^*}$. 
Problem~\ref{cond_expect_theta} is a minimization problem over a set of functions, which can be approximated by using risk minimization:
\begin{equation}
    \inf_\rho \sum_{i\in[n]}\mathbb{E}\left[\left(\chi(t_i,\hat{X}^{\theta^*}_{t_{i}};\rho)-Y^{\theta^*}_i\right)^2\right]
\end{equation}
 where $\chi(t,x;\rho)$ is deep neural network with parameter $\rho$; surrogate value function. Therefore,
the empirical risk minimization problem for approximation of the value function is given by
\begin{equation}\label{value_fnc:risk_minimization}
    \inf_\rho \sum_{i\in[n]}\sum_{m\in[M]}\left(\chi(t_i,\hat{X}^{m,\theta^*}_{t_{i}};\rho)-Y^{m,\theta^*}_i\right)^2
\end{equation}
where $\{\hat{X}^{m,\theta^*}:m\in[M]\}$ are sample trajectories and
\begin{equation}\label{eqn:Y_tmtheta*}
    Y^{m,\theta^*}_i:=\sum_{\hat{i}\in[n-i]}L\big(t_{i+\hat{i}},\hat{X}^{m,\theta^{*}}_{t_{i+\hat{i}}}, \phi(t_{i+\hat{i}},\hat{X}^{m,\theta^{*}}_{t_{i+\hat{i}}};\theta^{*})\big)\delta + g\big(\hat{X}_T^{m,\theta^{*}}\big)
\end{equation}
is the corresponding cost along the sample trajectories.
The surrogate for the \emph{value function} is given by a minimizer of \eqref{value_fnc:risk_minimization}, which we denote by $\chi(t,x;\rho^*)$.

%%%%%%%%%%%%%%%%%%%%%%%%%%%%%%%%%%%%%
\subsection{Cost analysis of Deep PGM}
\label{sec:cost}
In this section, we estimate the computational cost of running deep PGM, by counting the number of operations needed to perform gradient evaluation in the gradient descent algorithm.
We start by evaluating the number of operations required for performing back-propagation on $L_i^m(\theta):=L(t_{i},\hat{X}^{m,\theta}_{t_{i}},\pi(t_{i},\hat{X}^{m,\theta}_{t_{i}};\theta))$
in \eqref{prob:empirical_risk_minimization_theta}. By using chain rule, we obtain
\begin{equation}\label{L_theta}
\begin{split}
        \frac{\partial}{\partial\theta} L_i^m(\Upxi_i^m)=&(L_x(\Upxi_i^m)+L_{\pi}(\Upxi_i^m)\pi_x(t_{i},\hat{X}^{m,\theta}_{t_{i}};\theta))\frac{\partial}{\partial\theta}\hat{X}^{m,\theta}_{t_{i}}+L_\pi(\Upxi_i^m)\pi_\theta(t_{i},\hat{X}^{m,\theta}_{t_{i}};\theta)\\
        \frac{\partial}{\partial\theta} g(\hat{X}_T^{m,\theta}) = & \nabla g(\hat{X}_T^{m,\theta})\frac{\partial}{\partial\theta}\hat{X}^{m,\theta}_{T}
\end{split}
\end{equation}
where we set $\Upxi_i^m:=(t_{i},\hat{X}^{m,\theta}_{t_{i}},\pi(t_{i},\hat{X}^{m,\theta}_{t_{i}};\theta))$. Then, $\frac{\partial}{\partial\theta}\hat{X}^{m,\theta}_{t_{i}}$ is evaluated recursively by 
\begin{equation}\label{hatX_theta}
    \begin{split}
        \frac{\partial}{\partial\theta}\hat{X}^{m,\theta}_{t_{i+1}}=& \Big(I_d+(\mu_x(\Upxi_i^m)+\mu_\pi(\Upxi_i^m)\pi_x)\delta+(\sigma_x(\Upxi_i^m)+\sigma_\pi(\Upxi_i^m)\pi_x)\Delta W_{{t_{i+1}}}^{m}\Big)\frac{\partial}{\partial\theta}\hat{X}^{m,\theta}_{t_{i}}\\
        &+\mu_\pi(\Upxi_i^m) \pi_\theta(t_{i},\hat{X}^{m,\theta}_{t_{i}};\theta)\delta+\sigma_\pi(\Upxi_i^m)\pi_\theta(t_{i},\hat{X}^{m,\theta}_{t_{i}};\theta)\Delta W_{{t_{i+1}}}^{m}
    \end{split}
\end{equation}
The number of operations to evaluate the partial derivatives $\mu_x$, $\mu_\pi$, $\sigma_x$, $\sigma_\pi$, $\pi_x$, and $\pi_\theta$ via automatic differentiation is a constant. If we denote the number of operations for evaluation of $\frac{\partial}{\partial\theta}\hat{X}^{m,\theta}_{t_{i-1}}$ by $C_{i-1}$, then $C_i=C_{i-1}+\bar{C}$, where $\bar{C}$ is the number of operations to evaluate partial derivatives $L_x$, $L_{\pi}$, $\mu_x$, $\mu_{\pi}$, $\sigma_x$, $\sigma_{\pi}$, $\pi_x$, $\pi_\theta$, and tensor multiplications in \eqref{L_theta} and \eqref{hatX_theta}. Note that $\bar{C}$ depends on the architecture of the neural network, including the number of parameters. Therefore, the total number of operations for performing back-propagation along each sample trajectory is given by $C_n=n\bar{C}$. Recall that the total number of paths is denoted by $M$. Therefore, we require $nM\bar{C}$ operations to perform the back-propagation on the cost function. Note that the same number of operations as in $L_i^m(\theta)$ is required for the terminal cost $g(\hat{X}_T^{m,\theta^{*}})$. If \eqref{prob:empirical_risk_minimization_theta} is handled by a gradient descent algorithm, the number of operations is multiplied by the number of epochs.
\begin{remark}[Brute-force PGM]
    We refer to the implementation of the PGM discussed in this section as the brute-force PGM. The brute-force PGM serves as a benchmark to compare the efficiency of the deep hierarchical PGM discussed in the next section.
\end{remark}
%%%%%%%%%%%%%%%%%%%%%%%%%%%%%%%%%%%%%
%%%%%%%%%%%%%%%%%%%%%%%%%%%%%%%%%%%%%
\section{Deep hierarchical PGM}
 \label{sec:main}
The goal of this section is to solve the discretized approximation \eqref{prob:control_discretized} for large $n$, more efficiently than applying directly the brute-force deep PGM with $n$ time steps, described in \eqref{prob:risk_minimization_theta}-\eqref{X^theta}. The idea of deep hierarchical PGM is explained in the following steps:

\noindent\textbf{Step 1.} Assume that $n=N_1N_2$ and 
 consider the discrete problem \eqref{prob:control_discretized}-\eqref{Euler-Maruyama} with $n=N_1$, $\delta_1:= \frac{T}{N_1}$, and $T_{i}:=i\delta_1$. We call this problem the \emph{coarse discrete problem} and assume that there exists $\theta_1$ such that
 \begin{equation}\label{prob:coarse}
    \theta_1 \in\mathop{\text{\normalfont argmin}}_{\theta}\hat{\mathcal{J}}^{N_1}(\theta)
\end{equation}
where $\hat{\mathcal{J}}^{N_1}(\theta)$ is defined in \eqref{prob:empirical_risk_minimization_theta}. In other words, $\pi^1=\phi(t,x;\theta_1)$ provides an approximately optimal \emph{coarse policy} for \eqref{prob:empirical_risk_minimization_theta}. The \emph{coarse sample trajectories} under this policy are given by \eqref{X^theta}:
\begin{align}
   \hat{X}^{m,\theta_1}_{T_{i+1}} =\hat{X}^{m,\theta_1}_{T_{i}}&+\mu({T_{i}},\hat{X}^{m,\theta_1}_{T_{i}},\phi(T_{i},\hat{X}^{m,\theta_1}_{T_{i}};\theta_1))\delta_1\nonumber\\
   & +\sigma({T_{i}},\hat{X}^{m,\theta_1}_{T_{i}},\phi(T_{i},\hat{X}^{m,\theta_1}_{T_{i}};\theta_1))\Delta W_{{T_{i+1}}}^{m}, ~\hat{X}^{m,\theta_1}_0\sim\mathcal{D}_0\label{X^1}
\end{align}
where $\Delta W_{{T_{i+1}}}^{m}=W_{{T_{i+1}}}^m-W_{{T_{i}}}^m$ and $m\in[M_1]$. Recall that the surrogate for the \emph{coarse value function} is given by \eqref{value_fnc:risk_minimization}, which we denote by $\chi(t,x;\rho_1)$ where
\begin{equation}\label{defn:value_^1}
    \rho_1\in\mathop\text{argmin}\limits_\rho \sum_{i\in[N_1]}\sum_{m\in[M_1]}\left(\chi(T_i,\hat{X}^{m,\theta_1}_{T_{i}};\rho)-Y^{m,\theta_1}_i\right)^2
\end{equation}
and, the same as \eqref{eqn:Y_tmtheta*},
\[
Y^{m,\theta_1}_i:=\sum_{\hat{i}\in[N_1-i]}L\big(T_{i+\hat{i}},\hat{X}^{m,\theta_1}_{T_{i+\hat{i}}}, \phi(T_{i+\hat{i}},\hat{X}^{m,\theta_1}_{T_{i+\hat{i}}};\theta_1)\big)\delta + g\big(\hat{X}_T^{m,\theta_1}\big)
\]

\noindent\textbf{Step 2.} For each coarse interval $[T_i,T_{i+1}]$, $i\in[N_1]$, we measure two quantities: the change in the trajectory distribution and the change in the coarse value function across the interval, using the empirical training domains in Definition~\ref{defn:empirical_training_domain}. To measure the change in the trajectory distribution, we use the Hausdorff distance between the empirical training domain at time $T_i$ and $T_{i+1}$:
\begin{equation}\label{defn:hausdorff_score}
    d_H^i:=\max\left\{\sup_{x\in D_{\text{\normalfont tr}}(T_i)}\inf_{y\in D_{\text{\normalfont tr}}(T_{i+1})}|x-y|,\ \sup_{y\in D_{\text{\normalfont tr}}(T_{i+1})}\inf_{x\in D_{\text{\normalfont tr}}(T_i)}|x-y|\right\}
\end{equation}
between $D_{\text{\normalfont tr}}(T_i)$ and $D_{\text{\normalfont tr}}(T_{i+1})$, which measures how much the distribution of the coarse trajectory shifts across $[T_i,T_{i+1}]$. The measure of change in the value function is defined by the empirical mean squared difference of the surrogate for coarse value function $\chi(t,x;\rho_1)$ between $T_i$ and $T_{i+1}$:
\begin{equation}\label{defn:msd_score}
    \operatorname{MSD}^i:=\frac{1}{|S^i|}\sum_{x\in S^i}\big(\chi(T_i,x;\rho_1)-\chi(T_{i+1},x;\rho_1)\big)^2
\end{equation}
where $S^i$ is a finite grid of points evenly spaced in
\[
\operatorname{conv}\big(D_{\text{\normalfont tr}}(T_i)\big)\cap\operatorname{conv}\big(D_{\text{\normalfont tr}}(T_{i+1})\big),
\]
the overlap of the convex hulls of the two consecutive empirical training domains. Large values of $d_H^i$ and $\operatorname{MSD}^i$ indicate that the coarse time step $\delta_1$ is too coarse on $[T_i,T_{i+1}]$ to capture the local behavior of the optimal trajectory, and we use these scores to select the subset $\mathcal{I}\subseteq[N_1]$ of intervals to be refined in Step 3. 
\begin{remark}
    The measures of changes in the distribution of the coarse trajectory and the coarse value function across $[T_i,T_{i+1}]$ in this step are only used for the purpose of examples in Section~\ref{sec:numerics} and are not by any means the only ways to assess which intervals meet the criteria to be included in the next step.
\end{remark}

\noindent\textbf{Step 3.} Let $\mathcal{D}^1_{T_{i}}$ be the distribution of $\hat{X}^{m,\theta_1}_{T_{i}}$ and set $\delta_2=\sfrac{\delta_1}{N_2}=\sfrac{T}{N_1N_2}$, $t_{i,k}:=i\delta_1+k\delta_2$. For each $i\in\mathcal{I}$, we draw $M_2$ new samples from $\{\hat{X}^{m,\theta_1}_{T_i}:m\in[M_1]\}$ and use them with $M_2$ new sample paths of Brownian motion to simulate sample paths of $\hat{X}^{m,\eta}_{t_{i,k}}$ inside the interval $[T_i,T_{i+1}]$, i.e., $\{\hat{X}^{m,\eta}_{t_{i,k}}:i\in\mathcal{I},~k\in[N_2],~m\in[M_2]\}$, defined by
\begin{align}
   & \hat{X}^{m,\eta}_{t_{i,k+1}}=\hat{X}^{m,\eta}_{t_{i,k}}+\mu({t_{i,k}},\hat{X}^{m,\eta}_{t_{i,k}},\psi(t_{i,k},\hat{X}^{m,\eta}_{t_{i,k}};\eta))\delta_2\nonumber\\
   &\hspace{3cm} +\sigma({t_{i,k}},\hat{X}^{m,\eta}_{t_{i,k}},\psi(t_{i,k},\hat{X}^{m,\eta}_{t_{i,k}};\eta))\Delta W_{{t_{i,k+1}}}^m,~\hat{X}^{m,\eta}_{i,0}\sim\mathcal{D}^1_{T_i}\label{X^eta}
\end{align}
where $\Delta W_{{t_{i,k+1}}}^m=W_{t_{i,k+1}}^m-W_{t_{i,k}}^m$, and $\psi(t,x;\eta)$ is a deep neural network with parameter $\eta$ and possibly different architecture than $\phi$.
Then, we formulate the following empirical risk minimization problem to generalize the coarse optimal policy to the fine scale on the selected intervals $\mathcal{I}$:
\begin{equation}\label{prob:fine}
\begin{split}
    & \inf_{\eta}\sum_{i\in\mathcal{I}}\hat{\mathcal{J}}^{i,N_2}(\eta)\\
    &\hat{\mathcal{J}}^{i,N_2}(\eta):=\sum_{m\in[M_2]}\left[\sum_{k\in[N_2]}L(t_{i,k},\hat{X}^{m,\eta}_{t_{i,k}},\psi(t_{i,k},\hat{X}^{m,\eta}_{t_{i,k}};\eta))\delta_2 + \chi(T_{i+1},\hat{X}^{m,\eta}_{T_{i+1}};\rho_1)\right]
    \end{split}
\end{equation}
If $\mathcal{I}$ is evenly distributed in $[N_1]$, for instance $\mathcal{I}$ is the even numbers in $[N_1]$, then we expect to have a good interpolation property for $\psi(t,x;\eta^*)$ on $[T_{i},T_{i+1}]$ for $i\notin \mathcal{I}$. This study shows that we can leverage the different number of samples and different neural network architecture in each step to enhance the performance of the PGM.

%%%%%%%%%%%%%%%%%%%%%%%%%%%%%%
\subsection{\texorpdfstring{$K$}{}-fold implementation}
One can continue the deep hierarchical PGM beyond two steps to obtain finer time-discretization of the optimal control problem. 
Set $N_{k}=N\cdots N$ ($k$ times), for $k\in[K+1]$ and $T_{k,i}=i\delta_k$, $\delta_k=\sfrac{T}{N_k}$ with $i\in[N_k+1]$. Note that $T_{k,i}=T_{k+1,Ni}$.
Assume that, in stage $k$, one has approximated or evaluated the value functions $V^k(T_{k,i},x)$ for the discrete problem \eqref{prob:control_discretized} with $n=N_k$. Furthermore, assume that the distribution of the optimal trajectories of the discretized process $\hat{X}$ in \eqref{Euler-Maruyama} for step $k$ is known, i.e., $\hat{X}_{T_{k,i}}\sim \mathcal{D}^k_{i}$. Then, at the $k+1$st stage of the hierarchical method, we deal with the following empirical risk minimization problem:
\begin{align}
    &\inf_{\eta} \sum_{i\in\mathcal{I}_k}\hat{\mathcal{J}}^{k+1,i}(\eta)\\
   & \hat{X}^{m,\eta}_{T_{k+1,\ell+1}}\!\!\!\!\!=\hat{X}^{m,\eta}_{T_{k+1,\ell}}+\mu({T_{k+1,\ell}},\hat{X}^{m,\eta}_{T_{k+1,\ell}},\psi(T_{k+1,\ell},\hat{X}^{m,\eta}_{T_{k+1,\ell}};\eta))\delta_{k+1}\nonumber\\
   &\hspace{1.5cm} +\sigma({T_{k+1,\ell}},\hat{X}^{m,\eta}_{T_{k+1,\ell}},\psi(T_{k+1,\ell},\hat{X}^{m,\eta}_{T_{k+1,\ell}};\eta))\Delta W^m_{{T_{k+1,\ell+1}}},~\hat{X}^{m,\eta}_{T_{k+1,Ni}}\sim\mathcal{D}^k_{i}
\end{align}
where $i\in\mathcal{I}_k\subseteq [N_k]$, $\ell\in Ni+[N]$, $\Delta W^m_{{T_{k+1,\ell+1}}}= W^m_{{T_{k+1,\ell+1}}}- W^m_{{T_{k+1,\ell}}}$, $M_{k+1}$ is the number of i.i.d samples drawn from each distribution $\mathcal{D}^k_{i}$ and for Brownian motion sample paths, and
\begin{equation}
\begin{split}
     \hat{\mathcal{J}}^{k+1,i}(\eta):=&\sum_{m\in[M_{k+1}]}\bigg[\sum_{\ell\in Ni+[N]}L(T_{k+1,\ell},\hat{X}^{m,\eta}_{T_{k+1,\ell}},\psi(T_{k+1,\ell},\hat{X}^{m,\eta}_{T_{k+1,\ell}};\eta))\delta_{k+1}\\
     &\hspace{4cm} +\chi_k\Big(T_{k+1,N(i+1)},\hat{X}^{m,\eta}_{T_{k+1,N(i+1)}} ;\rho_k \Big)\bigg],
\end{split}
\end{equation}
where $\chi_k(t,x ;\rho_k)$ is the surrogate for the value function obtained in step $k$.
%%%%%%%%%%%%%%%%%%%%%%%%%%%%%%
\subsection{Computational cost of deep hierarchical PGM}
In Section~\ref{sec:cost}, we show that the cost of coarse PGM problem \eqref{prob:coarse} is $C_1N_1M_1$ operations, where $C_1$ depends on the architecture of the deep neural network $\phi_1(t,x;\theta)$, including the number of parameters.
To evaluate the surrogate for the value functions, the cost of \eqref{value_fnc:risk_minimization} is $C_{\text{\normalfont surr}}N_1M_1$, where $C_{\text{\normalfont surr}}$ depends on the architecture of the deep neural network $\chi_1(t,x;\rho)$, which makes the total computational cost for the coarse problem to be $(C_{\text{\normalfont surr}}+C_1)N_1M_1$.
With abuse of notation, we write the total computational cost of the first step as $c_1N_1M_1$, where $c_1:=C_1+C_{\text{\normalfont surr}}$ represents the computational cost imposed by the architecture of the deep neural network for policy $\phi_1(t,x;\theta)$ and the deep neural network for the surrogate $\chi_1(t,x;\rho)$.
For the refined problem in the next step, the cost can be similarly evaluated by $C_2N_1I_2N_2M_2$, where $C_2$ depends on the architecture of the deep neural network $\phi_2(t,x;\eta)$ for the refined step and $I_2=\sfrac{|\mathcal{I}|}{N_1}$. $N_1I_2$ is the number of selected intervals. Note that the number of sample trajectories for the refined problem is $I_2N_1M_2$.

If deep PGM is implemented directly with $n=N_1N_2$ and with $M$ samples, then the cost would be $CN_2N_1M$, where $C$ depends on the architecture of the neural network for the policy.
To obtain better efficiency with the deep hierarchical PGM, we require that $c_2N_1I_2N_2M_2+c_1N_1M_1<<CN_2N_1M$, or equivalently $C_2I_2M_2+c_1\sfrac{M_1}{N_2}<<CM$. If we conveniently ignore the difference between $c_k$ and $C_k$ and replace $C_k$ by $c_k$ and set $C=c_k$, we require $M_2+\sfrac{M_1}{N_2}=\sfrac{M}{\gamma}$ to improve the cost of the implementation in the hierarchical PGM by a factor of $\gamma$. For the $K$-fold hierarchical PGM we assume that the cost of back-propagation for both policy and surrogate training at step $k$ is given by $c_k$. Additionally, we denote by $M_k$, $I_k$, and $c_k$, respectively, the number of samples, the ratio of intervals included in the training over all intervals, and the cost of back-propagation at stage $k$. We note that $I_1=1$ always. We have the following result.

\begin{theorem}\label{thm:management}
For the $K$-fold hierarchical PGM with $N^K$ time steps to be $\gamma$ times more efficient than the brute-force PGM with the same number of time steps, it is sufficient to set $M_k$, and $c_k$, $k\in[K]+1$ such that    
\begin{equation}
    \begin{split}
        &a_K=\dfrac{1}{\gamma}-\dfrac{g_{K-1}}{N}, \text{and }~
    a_k=g_{k}-\dfrac{g_{k-1}}{N}~\text{ for }k\in[K]+1\\
    &\text{ with }~~ 0<\dfrac{g_{1}}{N^{K-1}}<\cdots<\dfrac{g_{K-1}}{N}<g_K=\dfrac{1}{\gamma}
    \end{split}
\end{equation}
where $a_k=\dfrac{c_{k}M_{k}I_{k}}{cM}$.
\end{theorem}
\begin{proof}
For $K$-fold hierarchical PGM, learning in stage $k+1$ includes $M_{k+1}I_{k+1}N_k$ samples, i.e., $I_{k+1}N_k$ is the number of intervals $[T_{k,i},T_{k,i+1}]$ used in training and $M_{k+1}$ is the sample size for each interval of the $I_{k+1}N_k$ intervals. Therefore, the computational cost for step $k+1$ is given by $c_{k+1}I_{k+1}M_{k+1}N_kN=c_{k+1}M_{k+1}I_{k+1}N_{k+1}$, where $c_{k+1}$ depends on the architecture of the neural network $\psi(t,x;\eta)$ at step $k+1$. Therefore, the total computational cost is given by $\sum_{k\in[K+1]}c_{k}M_{k}I_{k}N_{k}$ and the ratio of the cost relative to the cost of brute-force single stage PGM with $n=N_K$ and $M$ samples is
\begin{equation}
    \dfrac{\sum_{k\in[K]+1}c_{k}M_{k}I_{k}N_{k}}{cN_KM}=\!\!\!\!\!\sum_{k\in[K]+1}\!\!\!a_{k}\delta^{K-k},~\text{where}~ a_0=0, ~a_{k}=\dfrac{c_{k}M_{k}I_{k}}{cM}\text{for }k\ge1
\end{equation}
where $\delta=\delta_1=\sfrac{T}{N}$.
Note that $\sum_{k\in[K]+1}a_{k}\delta^{K-k}$ is the coefficient of $x^K$ in $g(x):=\sfrac{f(x)}{(1-\delta x)}$, where $f(x)=\sum_{k\in[K+1]}a_{k}x^{k}$.
In particular, we have $g(0)=0$ and $f(x)=(1-\delta x)g(x)$ must have positive coefficients. If we write $g(x)=\sum_{k=0}^\infty g_{k}x^{k}$, then $a_k>0$ and $g(0)=0$ imply
\begin{equation}\label{g g delta}
g_0=0,~g_{k+1}>\delta g_{k}~\text{ for }~k\in[K],~\text{ and }~g_{k}=\delta g_{k-1}~\text{ for }~k>K
\end{equation} 
Therefore, we must have 
\begin{equation}
    g(x)=\sum_{k\in[K]} g_{k}x^{k} +{g_Kx^K}\sum_{k\ge0}(\delta x)^k=\sum_{k\in[K]} g_{k}x^{k}+\dfrac{g_Kx^K}{1-\delta x}
\end{equation}
To improve the cost of $K$-fold hierarchical PGM $\gamma$ times, one needs to have $g_K= \sfrac{1}{\gamma}$. By \eqref{g g delta}, 
for any sequence of positive numbers $g_1,...,g_{K-1}$ such that \eqref{g g delta} is satisfied, the choice of $c_{k}M_{k}$ such that 
\begin{equation}
    \dfrac{c_{k}M_{k}I_{k}}{cM}=a_{k}=g_k-\dfrac{g_{k-1}}{N}
\end{equation}
yields $\gamma$-times more cost efficiency.
\end{proof}
\begin{eg}
    For $1$-fold PGM and $N=10$, one can choose $M_1=M$, $I_2=\frac{3}{10}$, $M_2=\frac{M}{2}$. Hence,
\begin{equation}
    \dfrac{c_{1}M_{1}}{cM}=a_1=1\quad\text{and}\quad\dfrac{c_{2}I_{2}M_{2}}{cM}=a_2=\dfrac{3}{20},
\end{equation}
which, by the recursion $g_k=a_k+\sfrac{g_{k-1}}{N}$, gives $g_1=1$ and $g_2=\dfrac14$. Therefore, using the brute-force with sample count of $M_1$ as benchmark, we gain $\gamma=\dfrac{1}{g_2}=4$ times in efficiency predicted by Theorem~\ref{thm:management}. In Section~\ref{sec:2-fold}, we see that implementation of this 2-fold PGM results in $4.2$ times empirical efficiency gain while keeping the same accuracy as the benchmark.
\end{eg}
\begin{eg}
    For $2$-fold PGM and $N=5$, one can choose $M_1=M$, $I_2=\frac25$, $M_2=\frac{M}{2}$, $I_3=\frac25$, $M_3=\frac{M}{2}$. Hence,
\begin{equation}
    \dfrac{c_{1}M_{1}}{cM}=a_1=1,\quad \dfrac{c_{2}I_{2}M_{2}}{cM}=a_2=\dfrac15,\quad\text{and}\quad\dfrac{c_{3}I_{3}M_{3}}{cM}=a_3=\dfrac15,
\end{equation}
which, by the recursion $g_k=a_k+\tfrac{g_{k-1}}{N}$, gives $g_1=1$, $g_2=\sfrac25$, and $g_3=\sfrac{7}{25}$. Therefore, using the brute-force with sample count of $M_1$ as benchmark, we gain $\gamma=\sfrac{1}{g_3}=\sfrac{25}{7}\approx3.57$ times in efficiency predicted by Theorem~\ref{thm:management}. In Section~\ref{sec:3-fold}, we see that implementation of this 3-fold PGM results in $5.2$ times empirical efficiency gain while keeping the same accuracy as the benchmark.
\end{eg}
  
%%%%%%%%%%%%%%%%%%%%%%%%%%%%%%
%%%%%%%%%%%%%%%%%%%%%%%%%%%%%%
%%%%%%%%%%%%%%%%%%%%%%%%%%%%%%
\section{Numerical experiment}
\label{sec:numerics}
All wall-clock timing results reported in this section were obtained on a single machine, a MacBook Air with an Apple M5 chip and 16 GB of memory. Unless noted otherwise, wall-clock statistics are computed over all completed independent runs of a given experiment; runs that failed to converge, or an extreme outlier runtime specifically identified in a table's footnote, are excluded from that table's statistics and reported separately. This is why the number of runs underlying each table differs.

\subsection{A linear-quadratic stochastic control problem}
\label{sec:lqsc}
In this section, we apply the proposed method on the standard continuous-time linear quadratic stochastic control (LQSC) problem below:
\begin{equation}\label{LQSC}
\begin{split}
    \inf_\pi~&\mathbb{E}\bigg[\int_0^T\big(aX_t^2+bX_t+A\pi_t^2+B\pi_t \big)dt+\alpha X_T^2+\beta X_T\bigg]\\
    &dX_{t}=(p X_t+q \pi_t)dt+\sigma dW_t
\end{split}
\end{equation}
The optimal policy for the LQSC \eqref{LQSC} can be given in closed form by $\pi^*_t=-\frac{B+q V_x(t,X^*_t)}{2A}$, where $X^*$ satisfies the SDE $dX^*_{t}=(p X^*_t+q \pi^*_t)dt+\sigma dW_t$ for the optimal trajectory and
$V(t,x)$ is the value function for \eqref{LQSC}, which is also given by a closed form $V(t,x)=f(t)x^2+h(t)x+k(t)$ with
\begin{equation}\label{Riccati}
         \begin{cases}
             0=f'+a+2pf-\frac{q^2}{A}f^2& f(T)=\alpha\\
             0=h'+b+ph-\frac{(B+qh)q}{A}f& h(T)=\beta\\
             0=k'+\sigma^2 f-\frac{1}{4A}(B+qh)^2& k(T)=0
         \end{cases}
\end{equation}
The system of ODEs in \eqref{Riccati} can be decoupled as the ODE for $f(t)$ is a Riccati equation, ODE for $h(t)$ is an integration away from $f(t)$, and ODE for $k(t)$ is an integration away from $f(t)$ and $h(t)$. More precisely, 
\[
\pi^*_t=-\tfrac{B+q(2f(t)X_t^*+h(t))}{2A}
\]
% \begin{remark}[Discretization bias in convex problems]
% In some problems, we may obtain $\hat{V}\le V$. For example, in the linear quadratic control problem with $\mu(x,u)=\alpha x+\beta u$ with $\alpha>0$ and $\beta<0$, $\sigma$ is constant, $L(x,u)=x^2+y^2$, and $g(x)=x^2+x$, the Jensen inequality implies that the discretized problem yields a smaller value function. In such cases, we comparing to the continuous-time closed-form solution is misleading. Therefore, we only use the closed-from solution to highlight the error of the brute-force method, while we compare the hierarchical with brute-force method with the same frequency.
% \end{remark}
%%%%%%%%%%%%%%%%%%%%%%%%%%%%%%%%%%%%%%%%%%%%%%%
%%%%%%%%%%%%%%%%%%%%%%%%%%%%%%%%%%%%%%%%%%%%%%%
%%%%%%%%%%%%%%%%%%%%%%%%%%%%%%%%%%%%%%%%%%%%%%%
\subsubsection{One-fold hierarchical experiment}
\label{sec:2-fold}
\begin{figure}[!htp]\centering
\includegraphics[width=\textwidth]{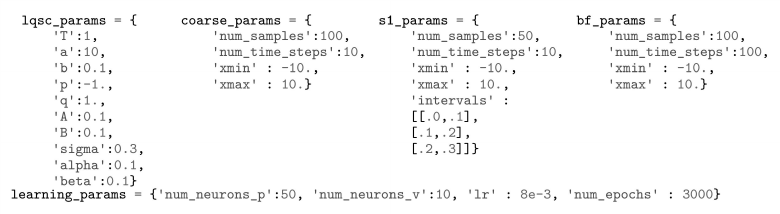}
\caption{\small The choice of parameters for the LQSC problem as well as the coarse (\texttt{lqsc}), the one-fold hierarchical method (\texttt{ms1}), and the brute-force (\texttt{bf}) method.}
    \label{params}
\end{figure}

The parameters of the LQSC problem are chosen as in Listing~\eqref{params} in \texttt{lqsc\_params}.

For the coarse step of this $1$-fold hierarchical PGM implementations, we take $N=10$, $M_1=100$, $X_0\sim\mathcal{D}_0=\text{Unif}(-10,10)$ as in Listing~\eqref{params} in \texttt{coarse\_params}.

\begin{remark}
    [Choice of parameters]
    We should emphasize that the problems of interest are the ones that a coarse discretization does not provide a relatively accurate solution for the continuous-time problem. The choice of parameters for the linear-quadratic stochastic control problem in our numerical experiment is such that the coarse discretization has a relatively large error and only when we refine the discretization, we obtain an accurate solution. To achieve that, we need functions $f(t)$ and $h(t)$ in \eqref{Riccati} to vary significantly in time. One way to achieve this variation is to choose the parameter $A$ relatively large, as we did in our case.
\end{remark}

For the coarse training, the surrogate for the policy is a deep neural network with two layers and $50$ neurons per layer. The number of epochs is set to 3,000. We consider the same parameters for training the surrogate for the coarse value function.

\begin{figure}[!htp]
\centering
    \includegraphics[width=1\linewidth]{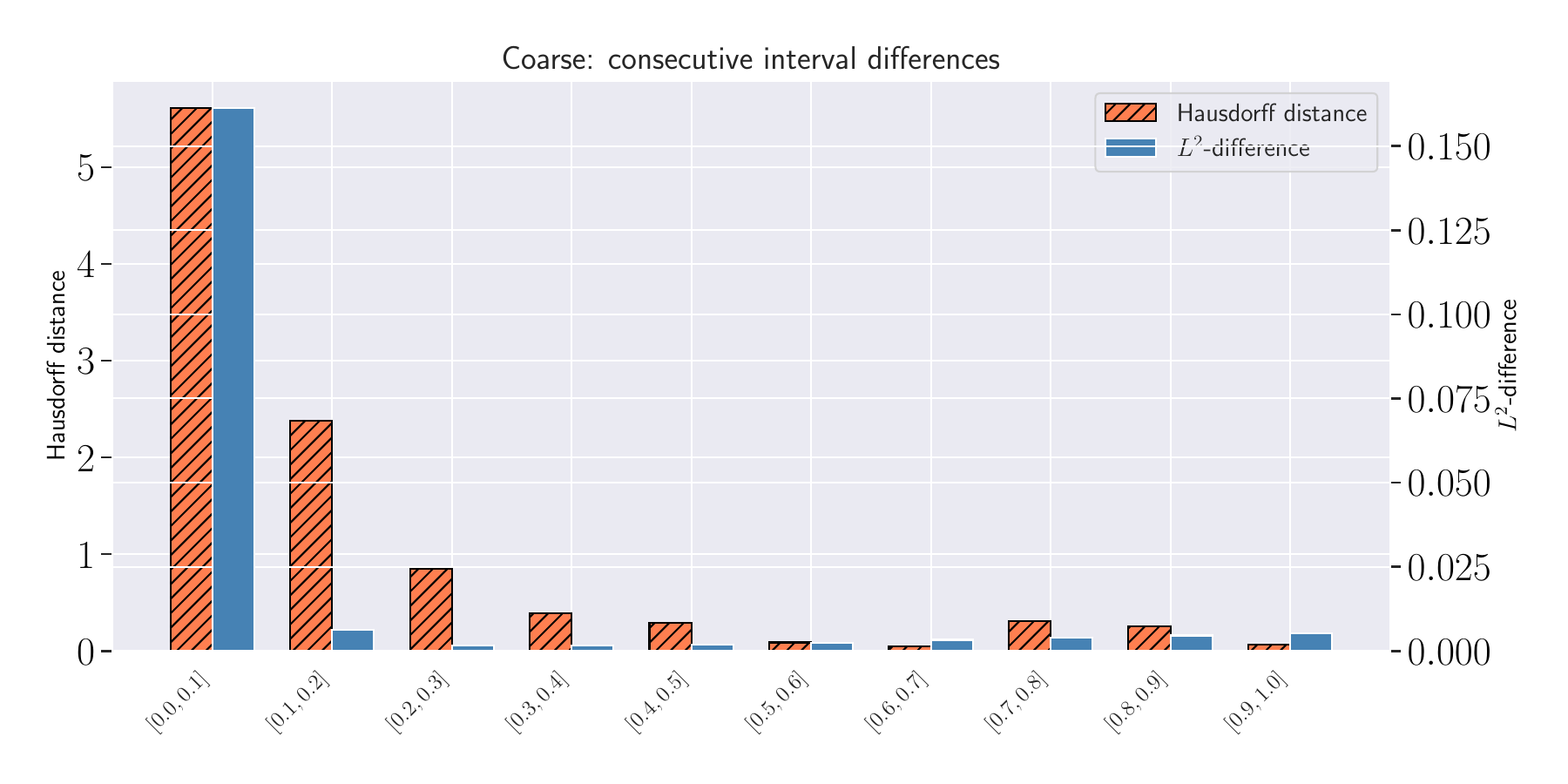}
    \caption{\small Right vertical axis: The mean squared difference between the surrogate value function $\chi(t_i,\cdot;\rho^*)$ and $\chi(t_{i+1},\cdot;\rho^*)$ on $D_{\text{\normalfont tr}}(t_i)\cap D_{\text{\normalfont tr}}(t_{i+1})$. Left vertical axis: The Hausdorff distance between $D_{\text{\normalfont tr}}(t_i)$ and $ D_{\text{\normalfont tr}}(t_{i+1})$.}
    \label{fig:interval_dist_ms1}
\end{figure}

In the fine step, we choose $M_2=50$ $N_2=10$ as in Listing~\eqref{params} in \texttt{s1\_params}, the surrogate for the fine policy is a deep neural network with two layers and $50$ neurons per layer, and the intervals are $[0,0.1]$, $[0.1,0.2]$, $[0.2,0.3]$. 
This choice is justified by the change between two consecutive time steps and the Hausdorff distance between the two consecutive empirical training domains shown in Figure~\ref{fig:interval_dist_ms1}.

We denote the coarse step by \texttt{lqsc}, the fine step by \texttt{ms1}, and the brute-force benchmark by \texttt{bf}.

\begin{figure}[!htp]
    \centering
    \includegraphics[width=0.75\linewidth]{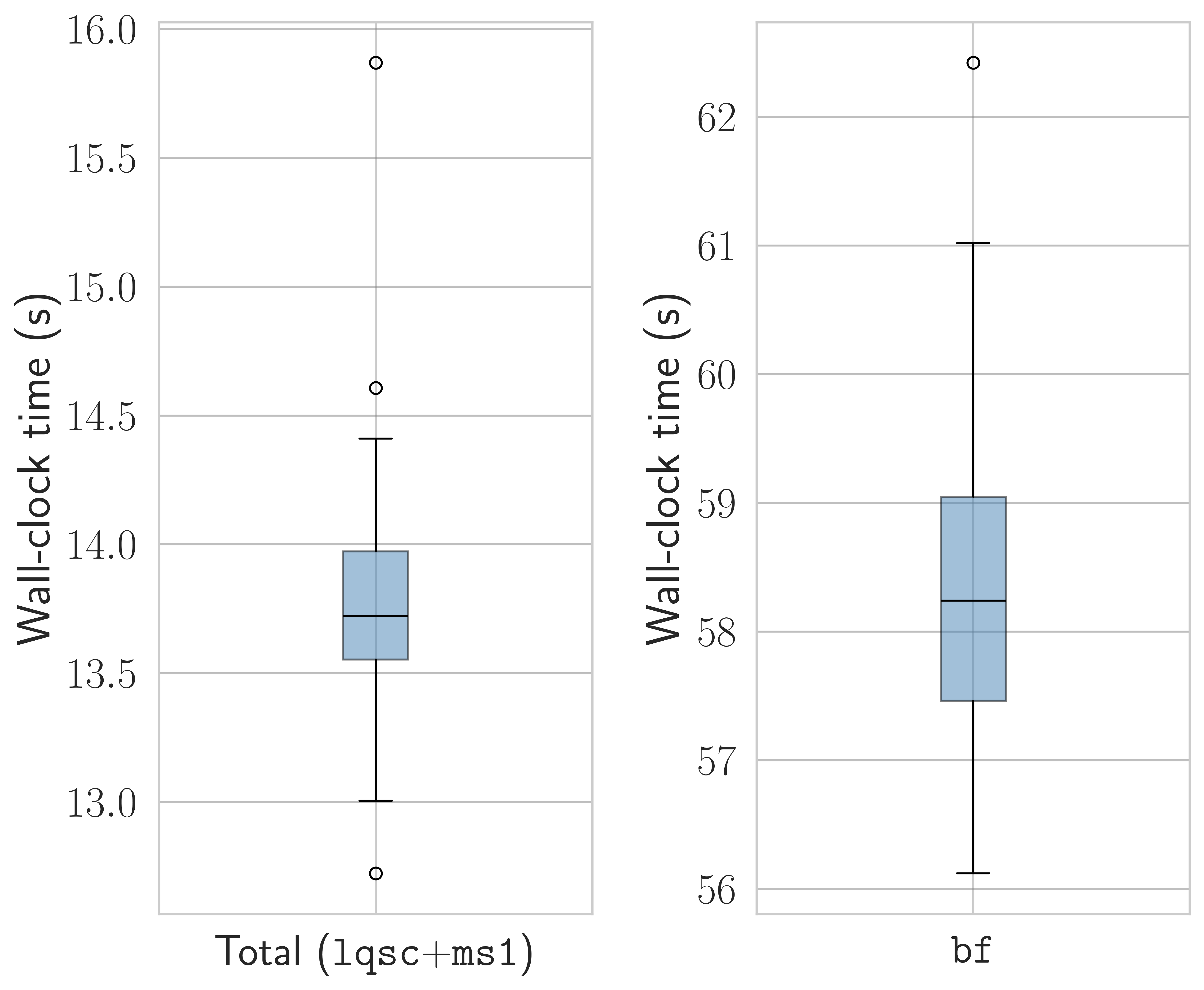}
    \caption{\small The wall-clock for one-fold PGM versus benchmark \texttt{bf} over 21 independent runs. There is an extreme outlier run for \texttt{bf} that took 401.23 s and is removed from the statistics.}
    \label{fig:lqc_2_wall_clock}
\end{figure}

\begin{figure}[!htp]
    \centering
    \includegraphics[width=0.95\linewidth]{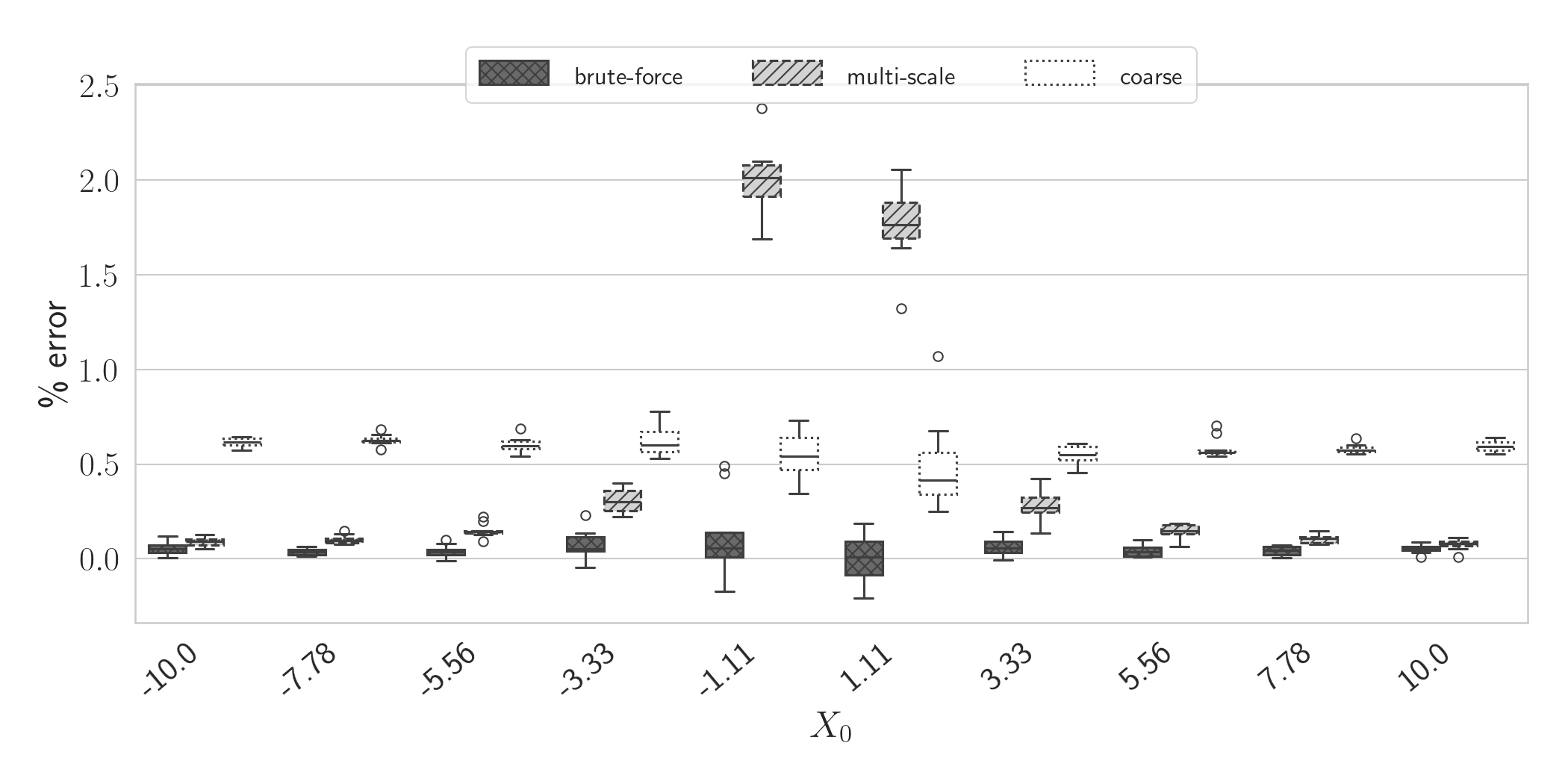}
    \caption{\small Comparison of the cost function between $1$-fold PGM and brute-force PGM with the closed-form solutions for $100$ time steps in 10 independent out-of-sample trajectories.}
    \label{fig:value0_comparison_bw_brute_multi}
\end{figure}

We measure the effectiveness of the fine surrogate policy by evaluating the average cost at $X_0\in \big\{\tfrac{i}{10}: i=-10,-9,...,10\big\}$ over 10 independent runs. As shown in Figure~\ref{fig:value0_comparison_bw_brute_multi}, the hierarchical method with one refinement in time discretization is significantly closer to the closed-form cost than the coarse loss and is comparable to the brute-force PGM in accuracy. The average running time of the hierarchical algorithm is around $4.2$ times faster. See Table~\ref{tab:training-times} and Figure~\ref{fig:lqc_2_wall_clock}.

\begin{table}[!htp]
\centering
\begin{tabular}{lcc}
\hline
Step & Avg. wall-clock time & Std dev \\
\hline
\texttt{lqsc.train()} ($N_1=10$) & 5.98 s & 0.35 s \\
\texttt{lqsc.value()} (not included in \texttt{lqsc.train()}) & 1.92 s & 0.36 s\\
\texttt{ms1.train()} (intervals $\{1,2,3\}$ \& $N_2=10$) & 5.88 s & 0.31 s \\
\hline
Total & 13.78 s & 0.64 s \\
\hline
\texttt{bf.train()}$^*$ & 58.43 s & 1.55 s \\
\hline
\end{tabular}
    \caption{\small Wall-clock training times for one-fold H-PGM in seconds over 21 independent runs. H-PGM is on average ${\sim}4.2\times$ faster than \texttt{bf}. $^*$ There is an extreme outlier run for \texttt{bf} that took 401.23 s and is removed from the statistics and the box plot.}
\label{tab:training-times}
\end{table}

% \begin{table}[!htp]
% \centering
% \begin{tabular}{llrrr}
% \hline
% \textbf{Method} & \textbf{Stage} & \textbf{Samples} & $N$ & \textbf{Time (s)} \\
% \hline
% \multirow{4}{*}{H-PGM}
%  & Coarse policy          & 100 & 10  & 5.73 \\
%  & Coarse value           & —   & —   & 0.92 \\
%  & Coarse value (retrain) & —   & —   & 1.45 \\
%  & ms1 policy             & 50  & 10  & 5.77 \\
%  & \textbf{Total}         &     &     & \textbf{13.87} \\
% \hline
% Benchmark & Brute-force  & 100 & 100 & 57.66 \\
% \hline
% \end{tabular}
% \caption{\small Wall-clock training times for LQC ($T=1$, 3{,}000 epochs each, one-fold H-PGM). H-PGM is ${\sim}4.5\times$ faster than brute-force. \texttt{ms1} refines the first 3 of the 10 coarse intervals.}
% \label{tab:wallclock_lqc_v7}
% \end{table}

%%%%%%%%%%%%%%%%%%%%%%%%%%%%%%%%%%%%%%%%%%

\subsubsection{Two-fold hierarchical experiment}
\label{sec:3-fold}
\begin{figure}[!htp]
\centering
    \includegraphics[width=0.95\linewidth]{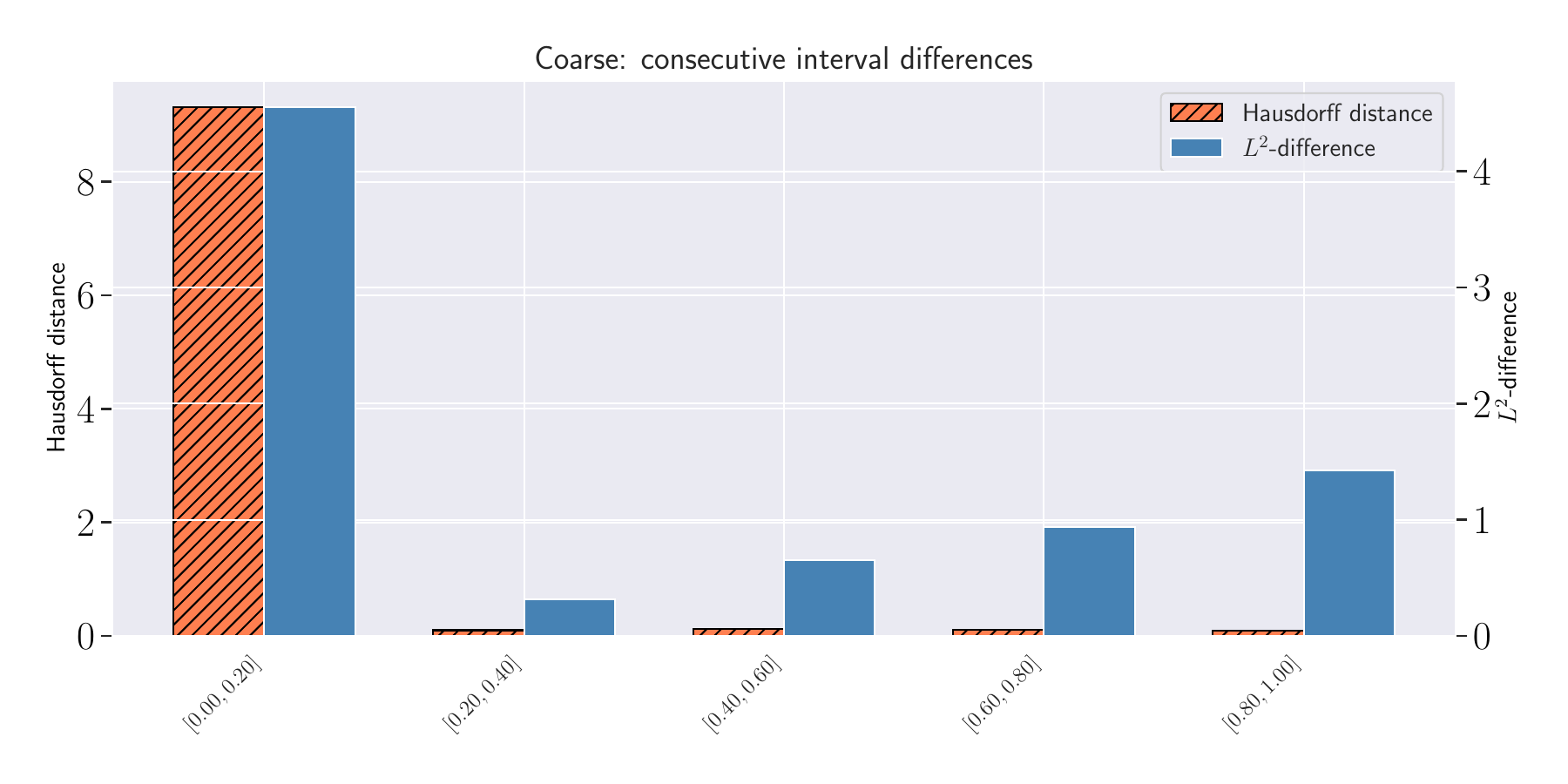}
    \caption{\small Right vertical axis: The mean squared difference between the surrogate value function $\chi(t_i,\cdot;\rho^*)$ and $\chi(t_{i+1},\cdot;\rho^*)$ on $D_{\text{\normalfont tr}}(t_i)\cap D_{\text{\normalfont tr}}(t_{i+1})$. Left vertical axis: The Hausdorff distance between $D_{\text{\normalfont tr}}(t_i)$ and $ D_{\text{\normalfont tr}}(t_{i+1})$.}
    \label{fig:L2_diff_consecutive_coarse_8}
\end{figure}

\begin{figure}[!htp]
\centering
    \includegraphics[width=0.95\linewidth]{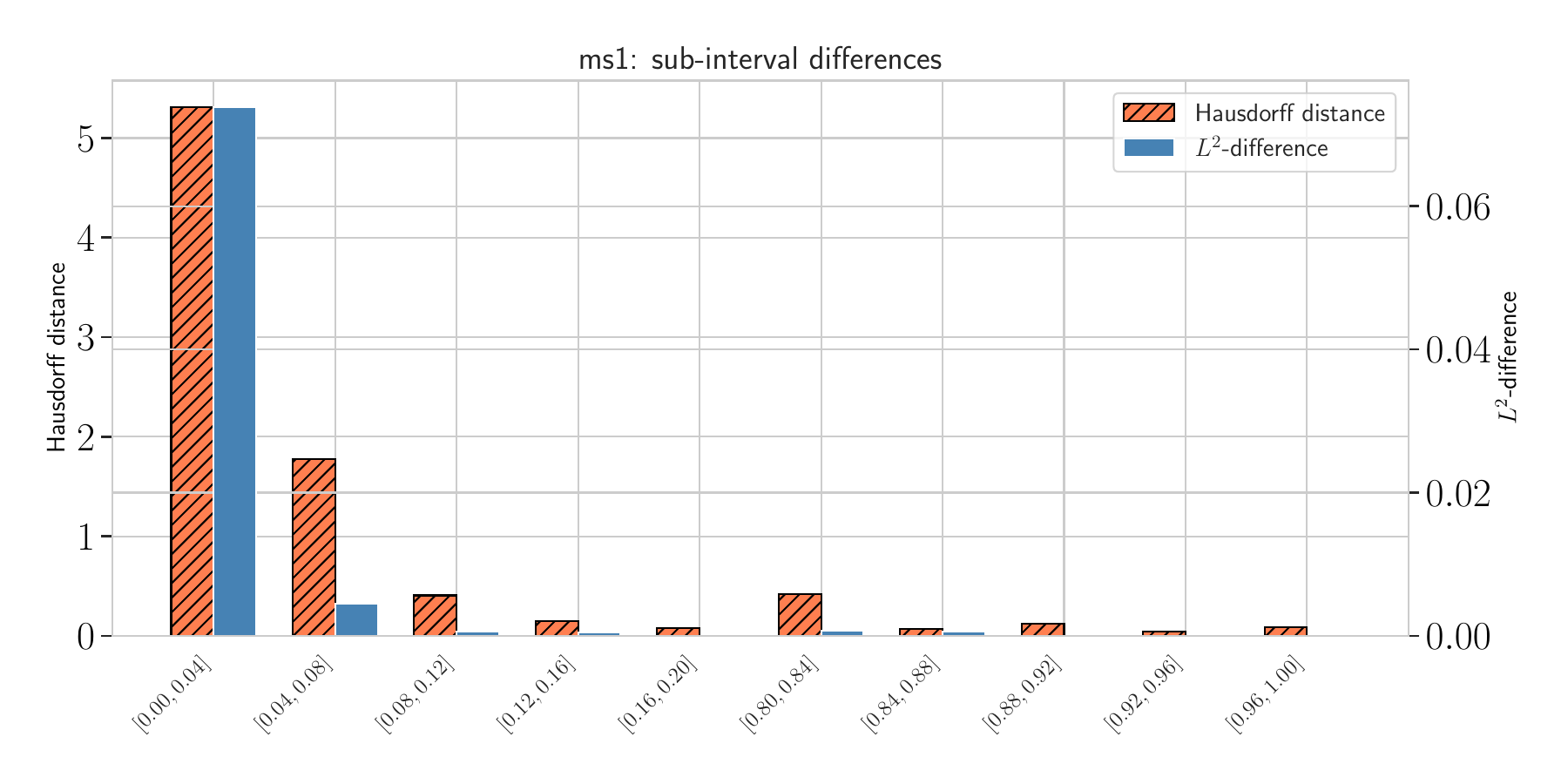}
    \caption{\small Right vertical axis: The mean squared difference between the fine surrogate value function $\chi(T_{i,k},\cdot;\rho^*)$ and $\chi(t_{i,k+1},\cdot;\rho^*)$ on $D_{\text{\normalfont tr}}(t_{i,k})\cap D_{\text{\normalfont tr}}(t_{i,k+1})$. Left vertical axis: The Hausdorff distance between $D_{\text{\normalfont tr}}(t_{i,k})$ and $ D_{\text{\normalfont tr}}(t_{i,k+1})$.}
    \label{fig:L2_diff_consecutive_ms1_8}
\end{figure}

\begin{figure}[!htp]\centering
\includegraphics[width=\textwidth]{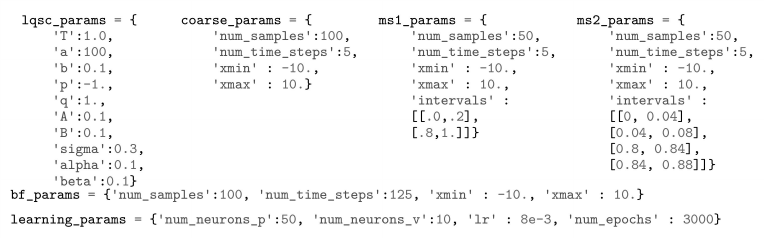}
\caption{\small The choice of parameters for the LQSC problem as well as the coarse (\texttt{lqsc}), the one-fold and two-fold hierarchical methods (\texttt{ms1} and \texttt{ms2}), and the brute-force (\texttt{bf}) method.}
    \label{params2}
\end{figure}

As seen in Listing~\eqref{params2}, for this experiment, we modify the coefficient $a$ from $10$ to $100$ to create more variation in time variable. We choose $N=5$, $M=100$, $M_1=100$, $M_2=50$, and $M_3=50$. Based on Figure~\ref{fig:L2_diff_consecutive_coarse_8}, the two coarse intervals $[0.0,0.2]$ and $[0.8,1.0]$ are chosen in the second step. For the training in the third step, we choose four intervals, two for each of the coarse intervals; $[0.0, 0.04]$ and $[0.04, 0.08]$ as subintervals of $[0.0,0.2]$, and $[0.8, 0.84]$ and $[0.84, 0.88]$ as subintervals of $[0.8, 1.0]$. These choices are made based on the variations observed in Figure~\ref{fig:L2_diff_consecutive_ms1_8}. The third step is denoted by \texttt{ms2}.

\begin{figure}[!htp]
    \centering
    \includegraphics[width=.95\linewidth]{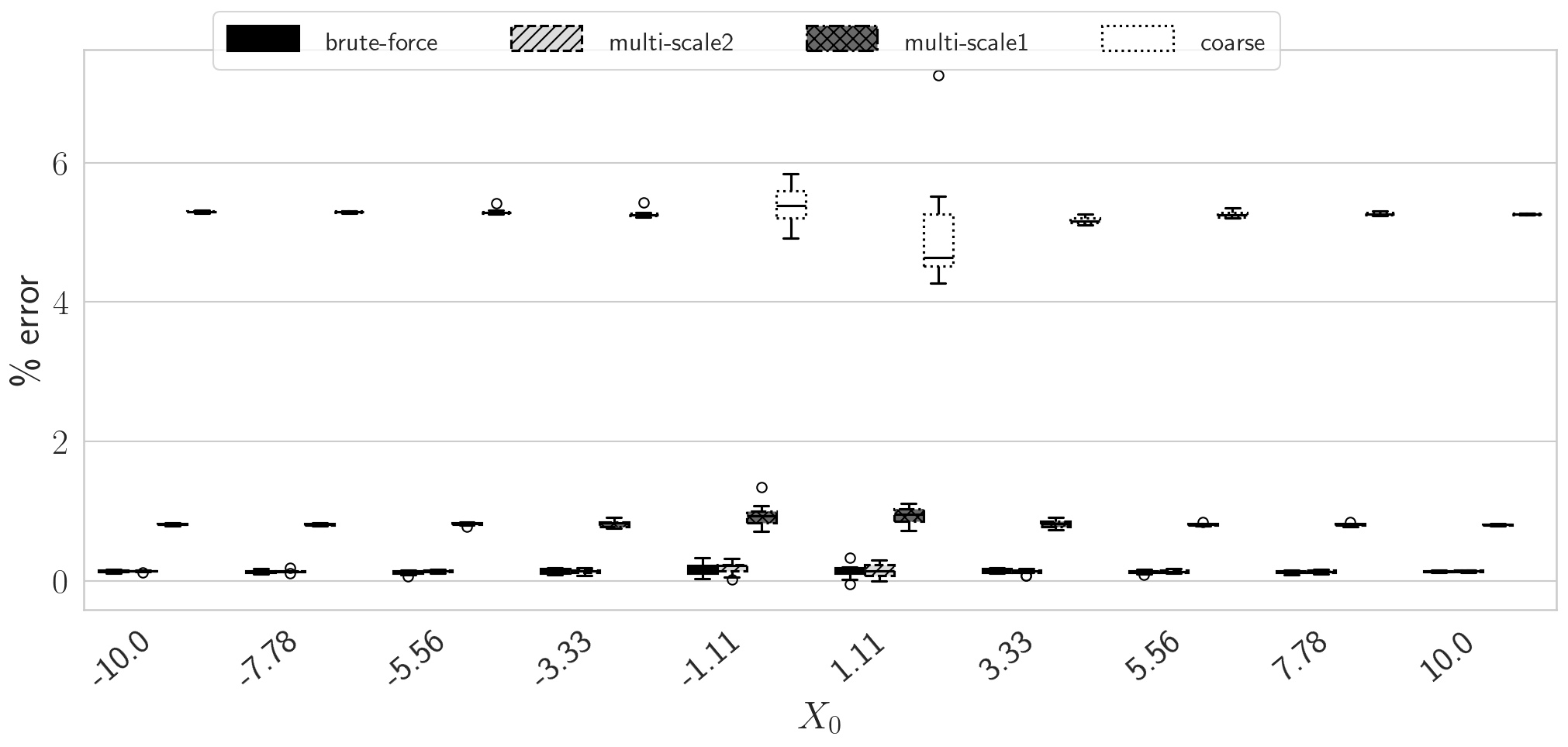}
    \caption{\small Error of the total expected cost of $2$-fold PGM and brute-force PGM for $125$ time steps in 10 out-of-sample trajectories at time $0$ over interval $[-10,10]$.}
    \label{fig:value2_comparison_bw_brute_multi}
\end{figure}

As seen in Figure~\ref{fig:value2_comparison_bw_brute_multi}, the error of the hierarchical method with two refinements in time discretization is comparable with the brute-force benchmark, while the efficiency gain is demonstrated in Table~\ref{tab:wallclock_lqc_v9}.

\begin{table}[!htp]
\centering
\begin{tabular}{lcc}
\hline
Step & Avg. wall-clock time & Std dev \\
\hline
\texttt{lqsc.train()} ($N_1=5$) & 3.2576 & 0.0905 \\
\texttt{lqsc.value()} (not included in \texttt{lqsc.train()}) & 3.7477 & 0.1237 \\
\texttt{ms1.train()} (intervals = first+last, $N_2=5$) & 3.0002 & 0.1093 \\
\texttt{ms1.value()} (not included in \texttt{ms1.train()}) & 0.4805 & 0.0362 \\
\texttt{ms2.train()} (4 sub-intervals, $N_3=5$) & 3.3345 & 0.0955 \\
\hline
Total & 13.8204 & 0.3395 \\
\hline
\texttt{bf.train()}$^*$ & 72.6379 & 1.4638 \\
\hline
\end{tabular}
\caption{\small Wall-clock training times for two-fold H-PGM in seconds over 20 independent runs. H-PGM is ${\sim}5.2\times$ faster than brute-force.}
\label{tab:wallclock_lqc_v9}
\end{table}

\begin{figure}[!htp]
    \centering
    \includegraphics[width=0.75\linewidth]{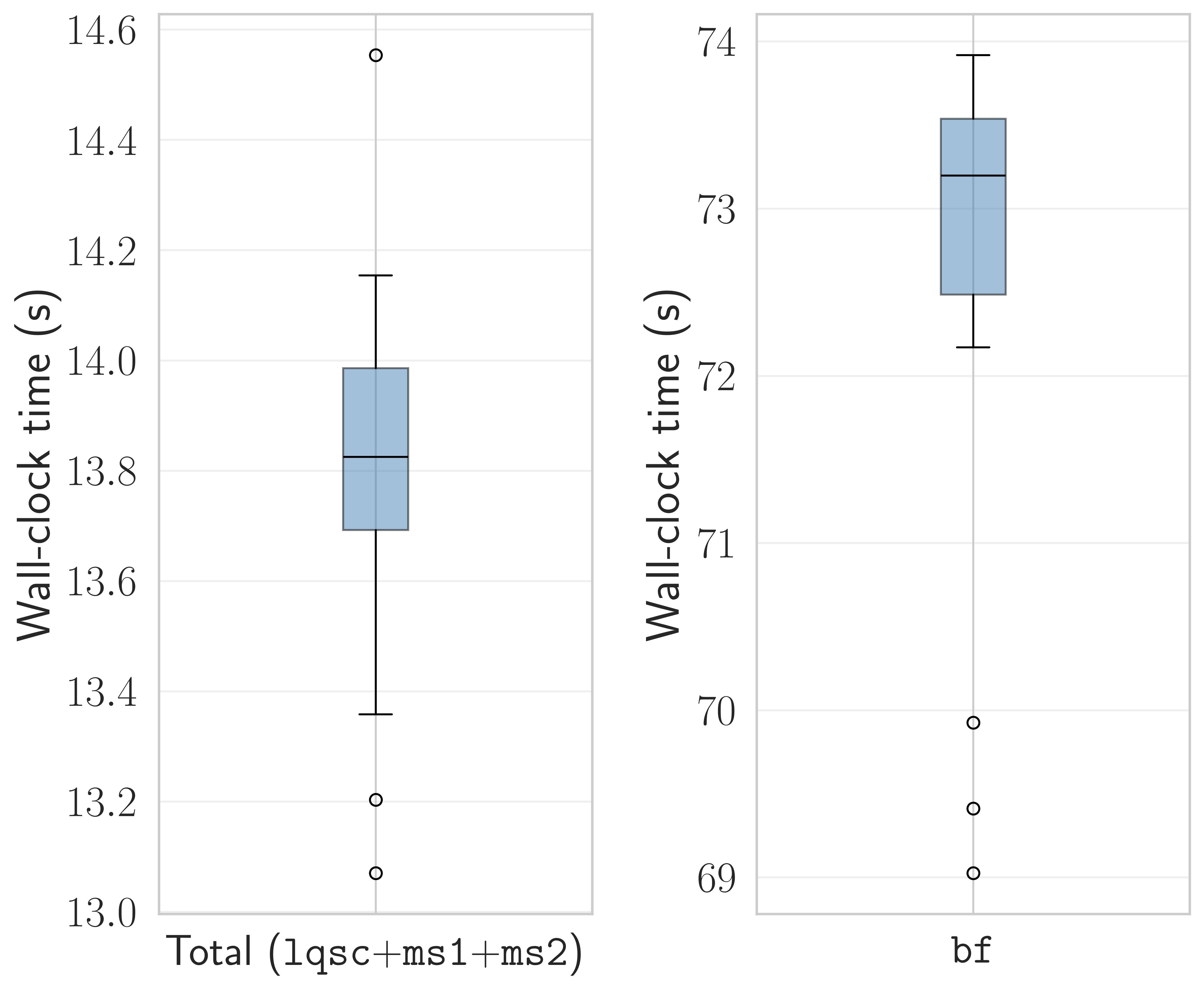}
    \caption{\small The wall-clock for two-fold PGM versus benchmark \texttt{bf} over 20 independent runs.}
    \label{fig:boxplot_3step}
\end{figure}

%%%%%%%%%%%%%%%%%%%%%%%%%%%%%%%%%%%%%%%%%%%%%%
%%%%%%%%%%%%%%%%%%%%%%%%%%%%%%%%%%%%%%%%%%%%%%

\subsection{Optimal execution under stochastic price impact and stochastic resilience}
\label{sec:oe}
The optimal execution problem arises when an entity is obligated to sell or buy a large quantity of an asset in a short period of time in the market. For example, when the composition of a market index, e.g., S\&P500, is revised, many account holders who use S\&P500 as a benchmark need to adjust their positions accordingly in a span of a few hours. This problem is discrete in nature. But, for practical purposes, it is required to be implemented in a higher frequency, which makes it closer to a continuous time problem. In many formulations, the continuous-time version of the problem is not interesting as the price impact component vanishes when the time step size goes to zero. 

We consider the linear price impact formulation in \cite{OW13}, which is generalized in many follow-up papers, \cite{AFS08,FSU14,FL17,FSU19}. The latter paper, which we use in this experiment, introduced stochastic price impact factor and resilience into the optimal execution problem. \cite{CLV24} also tackled a similar formulation which replaces the constraint on total liquidation by a balance penalty term and includes nonlinear price impact. For constant price impact factor and resilience and linear price impact, there is a closed-form solution, the well-known U-shaped strategy which makes two equally large trades at time $0$ and $N$ and equally small trades in between.  

For the purpose of testing hierarchical PGM, we consider the stochastic price impact and resilience model from \cite{FSU19}.
There are four state variables: the price impact $\kappa$, the price resilience $\rho$, temporary price impact $D$, and balance $R$.
We assume that the price impact $\kappa_t$ and resilience $\rho_t$ are given by OU processes
\[
d\kappa_t = \mu_1(\bar\kappa - \kappa_t)\,dt + \sigma_1\, dW_t^{1}
\]
\[
d\rho_t = \mu_2(\bar\rho - \rho_t)\,dt + \sigma_2\, dW_t^{2}
\]
where $W^1$ and $W^2$ are two independent Brownian motions, $\bar\kappa$ and $\bar\rho$ are the means of the stationary distributions of $\kappa_t$ and $\rho_t$, and $\mu_1, \mu_2, \sigma_1, \sigma_2$ are positive constants. The probability of $\kappa_t$ and $\rho_t$ taking negative values is non-zero. While negative values for $\rho_t$ do not cause degeneracy in the control problem, negative or zero values for $\kappa_t$ make the minimum cost $-\infty$ due to the concavity of the value function. Therefore, we introduce a lower bound for price impact. The price impact that enters into the cost function and temporary price impact dynamics is $\kappa_t\vee\underline{\kappa}$ for some $\underline{\kappa}>0$.
\begin{equation}
    \begin{split}
        \inf_{\xi_0,...,\xi_{N-1}} \mathbb E\bigg[&\sum_{i=0}^{N-1}\xi_i\Big(D_i + \frac{\kappa_{t_i}\vee\underline{\kappa}}{2}\,\text{\normalfont sgn}(\xi_i)|\xi_i|^\alpha\Big)\\
        &+ R_N \Big(D_N + \frac{\kappa_{T}\vee\underline{\kappa}}{2}\,\text{sgn}(R_N)|R_N|^\alpha\Big)\bigg]\\
        D_{i+1} = e^{-\rho_{t_i}\delta}\big(D_{i}+&(\kappa_{t_i}\vee\underline{\kappa}) \text{\normalfont sgn}(\xi_i)|\xi_i|^\alpha\big)\quad
    R_{i+1} = R_{i} -\xi_{i}, i=0,...,N-1
    \end{split}
\end{equation}
The terminal cost $g(D,R,\kappa,\rho)=DR+\tfrac{\kappa\vee\underline{\kappa}}{2}R^2$ indicates that we execute the whole remaining balance at time $N$. With $\alpha=1$, the price impact is linear: $D_{i+1} = e^{-\rho_{t_i}\delta}\big(D_{i}+(\kappa_{t_i}\vee\underline{\kappa})\xi_i\big)$.

\begin{figure}[!htp]
\centering
\includegraphics[width=0.99\linewidth]{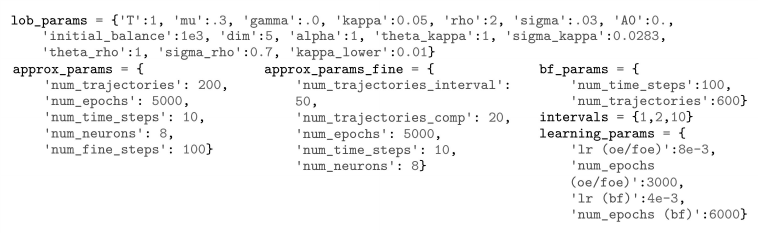}
    \caption{\small The parameter choice for the optimal execution numerical experiment.}
    \label{fig:lob_params}
\end{figure}

%%%%%%%%%%%%%%%%%%%%%%%%%%%%%%%%%%%%%%%%%%
\subsubsection{Degeneracy and challenge of training the value function surrogate}
\label{sec:degeneracy}
\begin{figure}[!htp]
\centering
    \includegraphics[width=0.98\linewidth]{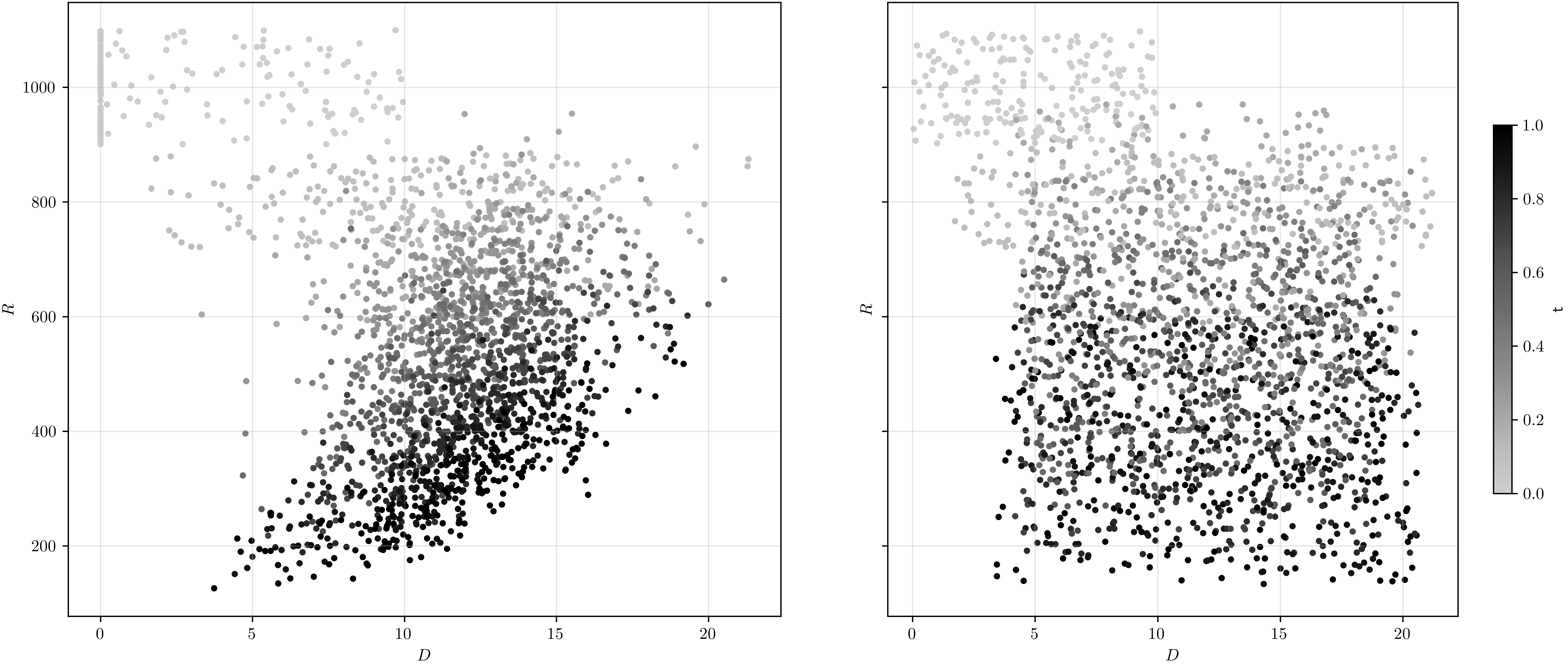}
    \caption{\small Left: points on the optimal coarse trajectories. To make it possible to train the value function and compare it at two consecutive time steps, we extend the training trajectories to larger rectangles at each time step. Different time steps are distinguished by different shades ranging from light to dark.}
    \label{fig:fat_data}
\end{figure}

One of the challenges of implementing hierarchical PGM on this problem is the degeneracy of variables $D$ and $R$. When training the coarse model, the optimal trajectories make a domain for the value function at each time step, which is a perturbed 3-dimensional manifold of coarse trajectories in the 4-dimensional state space $(D,R,\kappa,\rho)$. In left plot in Figure~\ref{fig:fat_data}, we see the points on the coarse trajectories distinguished in time step by different shades. Darker shades are the earlier time steps. The empirical training domains for different time steps have small overlap, and the domains become narrower for time steps near the terminal time. Specifically, there is a large correlation between $D$ and $R$. This causes the surrogate for the value function to quickly lose accuracy away from the perturbed 3-dimensional manifold. Even when we choose $(D_0,R_0)$ from a rectangle, at later time steps the collinearity emerges. Since the coarse value function is used in the fine-scale training, the fine training becomes significantly inaccurate as soon as a point on one of the fine sample trajectories deviates from the perturbed 3-dimensional manifold. It renders the naïve hierarchical PGM unstable. 

The remedy we propose is to add a new step in the implementation after the coarse training. We expand our training data at each time step to a rectangle which includes the perturbed 3-dimensional manifold. See right plot in Figure~\ref{fig:fat_data}. The new training data points do not suffer from collinearity. We run new coarse trajectories from each new data point according to the coarse optimal policy and evaluate the trajectory cost. Finally, we train a surrogate value function on the new training data points and their associated cost.

\begin{figure}[!htp]
\centering
    \includegraphics[width=0.95\linewidth]{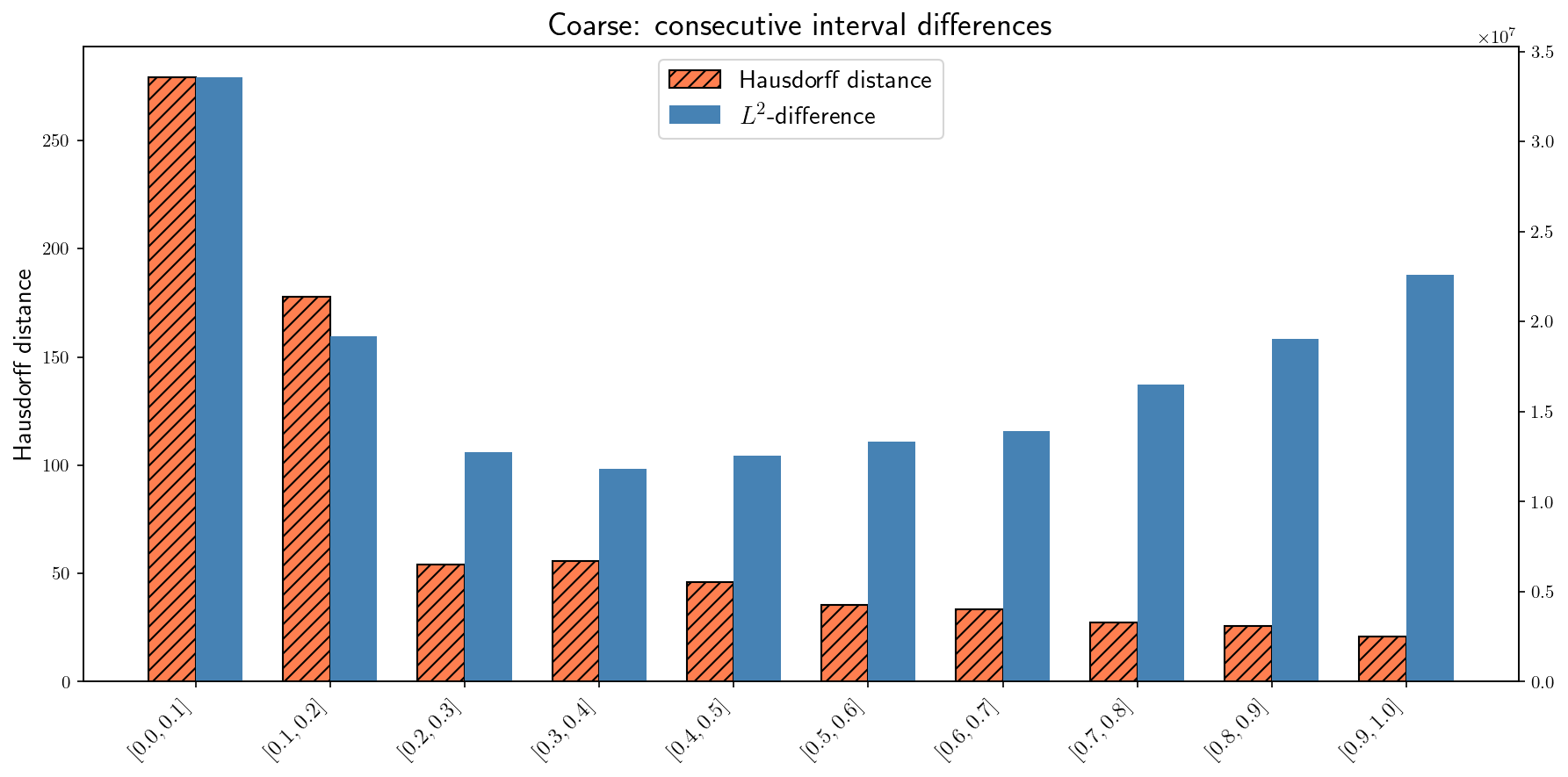}
    \caption{\small The comparative mean squared difference of value function at two consecutive time steps and Hausdorff distance between the consecutive empirical training domains.}
    \label{fig:L2_Hausdorff_consecutive_coarse_oe}
\end{figure}

As in Figures~\ref{fig:interval_dist_ms1} and \ref{fig:L2_diff_consecutive_coarse_8} for the LQSC experiments, Figure~\ref{fig:L2_Hausdorff_consecutive_coarse_oe} shows the Hausdorff distance and $L^2$-difference of the coarse surrogate value function between consecutive coarse intervals. The intervals $[0,0.1]$ and $[0.1,0.2]$ have the largest Hausdorff distance and $L^2$-difference and are selected for fine-scale refinement, together with $[0.9,1.0]$, where the $L^2$-difference rises again toward the terminal time. Since these three intervals concentrate near the two ends of the time horizon, we also include $[0.4,0.5]$, roughly midway between them, so that the fine-scale policy network generalizes more accurately across the interior of $[0,T]$.

\subsubsection{Total-liquidation constraint in the fine time scale}
Theoretically, the dynamic programming principle implies that the fine training must lead to total liquidation if the value function is accurately evaluated. For obvious reasons, this is not guaranteed in the hierarchical method. It is possible that a fine pretrained trajectory lies outside the domain of the value function. Apart from causing significant instability, this does not allow for approximate complete liquidation at the terminal time. To reduce the possibility of instability as well as enforcing approximate total liquidation, we include a small number of fine trajectories from time $0$ to the terminal time $T$ in the loss function for the fine training. 
More precisely, the fine minimization problem in \eqref{prob:fine} is now
\[
\inf_{\eta}\sum_{i\in\mathcal{I}}\mathcal{J}^{i,N_2}(\eta) + 
    \hat{\mathcal{J}}^{n}(\eta),~~~n=N_1N_2
\]
where $\hat{\mathcal{J}}^{n}(\eta)$ is given by \eqref{prob:empirical_risk_minimization_theta} with $M=M_0$, where $M_0$ is the number of fine trajectories from time $0$ to the terminal time $T$. The higher $M_0$ is, the better stability is gained, however, the efficiency is reduced. Without adding a collection of fine trajectories, we observed instability in half of the runs.

\subsubsection{The result of numerical implementation}

In the numerical implementations, we choose the parameters as in Figure~\ref{fig:lob_params}. The number of epochs listed in Figure~\ref{fig:lob_params} is the maximum number of epochs. In this experiment, we added a new feature that all trainings stop as soon as the relative change in the loss between two consecutive epochs is less than \texttt{1e-9}. In this problem, the value of loss is strictly positive and away from zero, which makes the evaluation of relative change in the loss between two epochs possible.

We run the hierarchical PGM $22$ times to record the results. Two of the $22$ runs were unstable. See Remark~\ref{rem:degeneracy2/20} for more details. \texttt{oe.train()} represents the coarse step. \texttt{oe.value()} represents evaluation of the coarse trajectories, expansion of coarse domain of the value function to combat degeneracy, and fitting of the value function. \texttt{foe.train()} represents the fine scale training result and \texttt{oe100} is the brute-force benchmark. The average total cost, average remaining balance, and the std of remaining balance are shown in Table~\ref{tab:cost_balance_stat}. As seen, the hierarchical PGM \texttt{foe} is comparable to the brute-force benchmark \texttt{oe100}. 

The relative-difference statistics between \texttt{foe} and \texttt{oe100} in the 20 runs are 
\[
\text{Mean rel. diff } \frac{\text{foe} - \text{oe100}}{|\text{oe100}|} = -0.2089\%, \quad \text{std} = 1.6716\%
\]
\[
\text{Min/Max rel. diff} = -3.6906\% \;/\; 2.8151\%
\]

\begin{table}[!htp]
\centering
\begin{tabular}{lcc}
\hline
 & \texttt{foe} & \texttt{oe100} \\
\hline
Mean total cost & 13299.33 & 13323.00 \\
Mean final $R$ & 285.49 & 273.49 \\
Std final $R$ & 49.10 & 50.63 \\
\hline
\end{tabular}
\caption{\small Over 20 independent runs: mean total cost and final remaining inventory $R$ for \texttt{foe} vs \texttt{oe100}.}
\label{tab:cost_balance_stat}
\end{table}

\begin{figure}[!htp]
    \centering
    \includegraphics[width=0.5\linewidth]{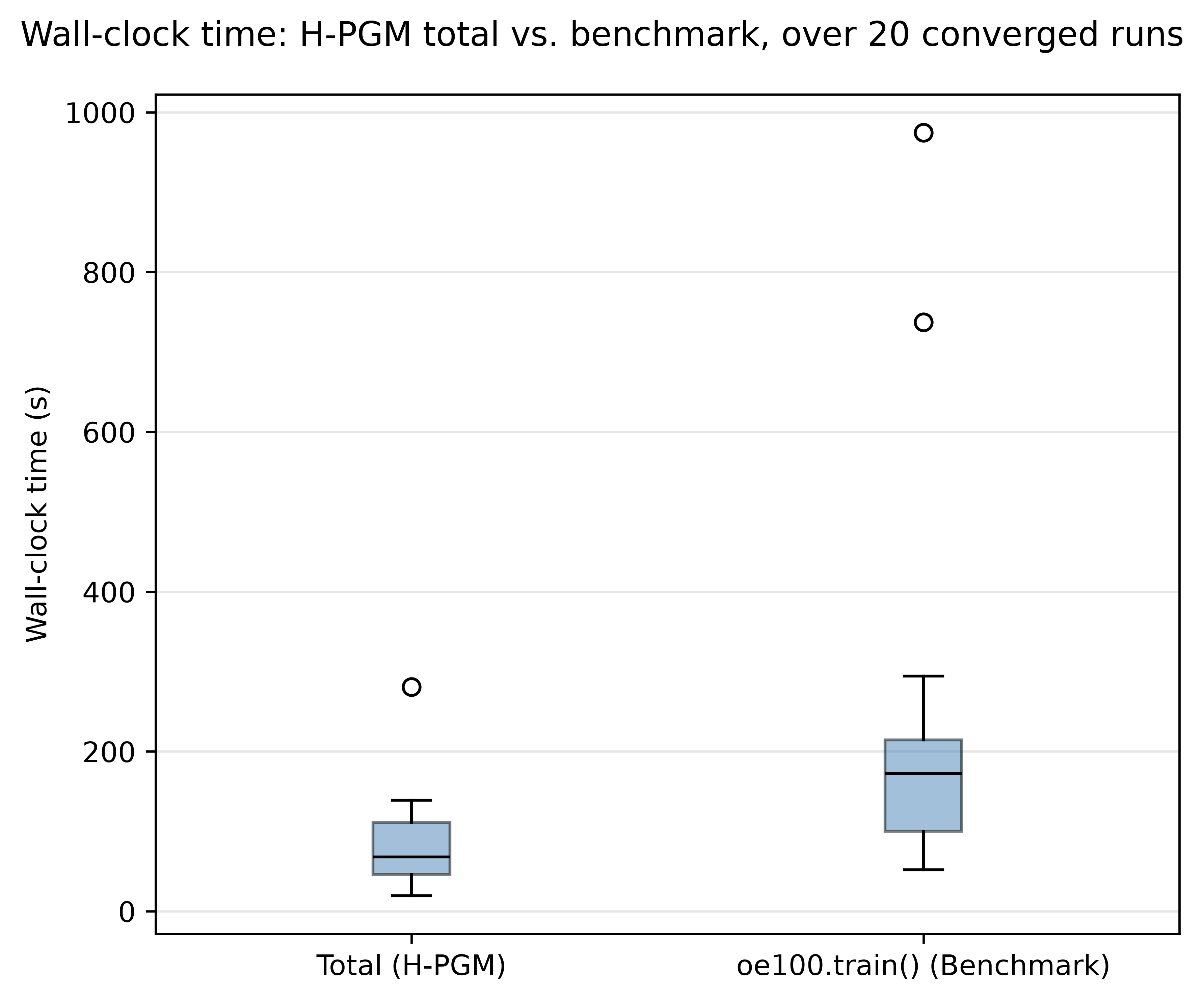}
    \caption{\small The boxplot of the wall-clock time across 20 runs.}
    \label{fig:boxplot}
\end{figure}

Table~\ref{tab:wallclock_lob_v7} shows the wall-clock time of running the hierarchical PGM versus the brute-force benchmark PGM. The box plot of the wall-clock time across 20 runs is shown in Figure~\ref{fig:boxplot}.

\begin{table}[!htp]
\centering
\begin{tabular}{lcc}
\hline
Step & Avg. wall-clock time & Std dev \\
\hline
\texttt{oe.train()} ($N_1=10$)& 15.90 s & 9.36 s \\
\texttt{oe.value()} (included in \texttt{oe.train()}) & 2.43 s & 1.01 s \\
\texttt{foe.train()} (intervals $\{1,2,5,10\}$ \& $N_2=10$) & 66.67 s & 56.30 s \\
\hline
Total& 82.57 s & 58.40 s \\
\hline
\texttt{oe100.train()} ($N=100$) & 224.19 s & 229.76 s \\
\hline
\end{tabular}
\caption{\small Mean and standard deviation of wall-clock time over 20 converged runs.}
\label{tab:wallclock_lob_v7}
\end{table}

It is quite possible that in some of the runs, the benchmark \texttt{oe100} stops earlier than \texttt{foe}. In fact, this happened in four of the 20 converging runs as shown in Table~\ref{tab:100>10x10}.
\begin{table}[!htp]
    \centering
    \begin{tabular}{cc}
\hline
Total H-PGM (s) & \texttt{oe100} (s) \\
\hline
96.87 & 75.41 \\
69.45 & 66.45 \\
67.38 & 51.75 \\
138.84 & 108.26 \\
\hline
\end{tabular}
    \caption{\small Four runs where the benchmark beat the hierarchical PGM in runtime.}
    \label{tab:100>10x10}
\end{table}

\begin{remark}[On training stability]
\label{rem:degeneracy2/20}
In a 22-run replication study, the hierarchical PGM
training pipeline (coarse \texttt{oe}/\texttt{oe100} and fine-scale \texttt{foe} policies) converged
to a stable solution in 20 of 22 runs. In the remaining 2 runs, one of the two policy networks diverged to a degenerate solution during training: in one case it was the
benchmark policy, in the other the fine-scale policy.
\end{remark}
%%%%%%%%%%%%%%%%%%%%%%%%%%%%%%%%%%%%%%%%%%%%%%
%%%%%%%%%%%%%%%%%%%%%%%%%%%%%%%%%%%%%%%%%%%%%%
\textbf{Data Availability Statement.} All the codes related to the numerical experiment performed in this study are available at \url{https://github.com/arashfahim/multiscale-PGM-for-stochastic-control}. The Jupyter notebooks for the $1$-fold and $2$-fold PGM are \url{LQC/Two-steps/LQC_Updated_v7.ipynb} and \url{LQC/Three-steps/LQC_Updated_v8.ipynb} and the optimal execution experiment is stored in \url{LOB/MS_PGM_v7.ipynb}.
% Several other experiments for different parameter sets are also available on this repository alongside experiment on other problems such as optimal execution under linear, nonlinear, and stochastic price impact.
%%%%%%%%%%%%%%%%%%%%%%%%%%%%%%%%%%%%%%%%%%%%%%
%%%%%%%%%%%%%%%%%%%%%%%%%%%%%%%%%%%%%%%%%%%%%%
\section{Conclusion}
In this paper, we proposed a method that improves the efficiency of the policy gradient method for optimal control problems via a hierarchical implementation. Our method is useful for problems where the coarse discretization does not provide a sufficiently accurate approximation and finer time discretization is required, specifically, if the cost varies significantly in time. We tested our method on a linear-quadratic stochastic control problem, which has a closed-form solution, and compared the results of one-fold and two-fold implementations. We also tested our method on an optimal execution problem with stochastic price impact and stochastic resilience. For the optimal execution problem, there are two challenges, degeneracy caused by the deterministic variables and total-liquidation constraint, which are tackled by introducing new sub-steps in the hierarchical PGM. Additionally, Theorem~\ref{thm:management} gives a theoretical estimate of the computational cost of the hierarchical algorithm as a function of the number of refinement folds, the fraction of intervals refined at each fold, and the number of samples used at each fold. The theorem provides a sufficient condition on these design choices that guarantees the hierarchical method achieves any prescribed speed-up factor $\gamma$ relative to the brute-force benchmark.

\section*{Acknowledgments}
The authors used AI as a writing and technical assistant in preparing this manuscript. Its use included: proofreading for grammar and typographical consistency; identifying and helping correct several errors and inconsistencies in the mathematical exposition (including equation labeling, index notation, and cross-references); and LaTeX formatting and typesetting assistance, including converting the manuscript to the SIAM LaTeX template. All AI-assisted suggestions were reviewed, verified, and, where mathematical, independently checked by the authors. The authors assume responsibility for all content of this paper. The first version of codes in the numerical experiments are created by the authors and later modified by the help of AI.

 \appendix
 \section{Proof of Theorem~\ref{thm:strong_error} and Proposition~\ref{prop:epsilon_optimal}}
 \label{sec:error}
 In this appendix, we provide a proof of the error bound for approximation of the continuous-time control problem with the discretized version stated in Theorem~\ref{thm:strong_error}. Although it is a standard result, we provide it for the sake of completeness. To establish the result, we use the estimate obtained by \cite[Theorem~2.1]{JPR19} for the restriction to piecewise-constant controls. We first assume the following:
 \begin{enumerate}[label=\bfseries A.\arabic*.]
     \item The set $\mathbb{U}\subset\mathbb{R}^m$ of values taken by admissible controls is compact.
     \item $L(t,x,u)$, $\mu(t,x,u)$, $\sigma(t,x,u)$, and $g(x)$ are bounded, Lipschitz in $x$, and uniformly in $(t,u)\in[0,T]\times\mathbb{U}$, and $L(t,x,u)$ is $\sfrac{1}{2}$-Hölder continuous in $t$ uniformly in $(x,u)$.
 \end{enumerate}
 \begin{remark}[Compactness and boundedness]
 Assumptions \textbf{A.1} and \textbf{A.2} are stated in the form required to invoke \cite[Theorem~2.1]{JPR19} below, which assumes a compact control-value set and bounded coefficients. Neither holds verbatim for the linear-quadratic problems of Section~\ref{sec:numerics}, where $u$ is unconstrained and $\mu,\sigma,L$ have linear or quadratic growth in $x$. As is standard in this literature, both can be recovered for the purposes of this proof without affecting the stated rate: $u$ is restricted to a compact set by saturation, and $\mu,\sigma,L,g$ are truncated outside a ball $\{|x|\le R\}$; a priori moment bounds on the optimal trajectory, available in closed form for linear-quadratic dynamics via the associated Riccati equation, then show that the truncation error is negligible relative to $\delta$ once $R=R(\delta)\to\infty$ at a suitable rate. We omit the details of this standard localization argument.
 \end{remark}
 \begin{proof}[Proof of Theorem~\ref{thm:strong_error}]
 We remind the definition of $\Uppi_\delta$ from Section~\ref{sec:intro}. For $\pi\in\Uppi_\delta$, let $X^\pi$ solve \eqref{prob:control} under $\pi$, and let $\hat{X}^\pi$ be its Euler-Maruyama discretization \eqref{Euler-Maruyama}, driven by the same Brownian path and the same control values $\{\pi_{t_j}:j\in[N]\}$. Define
 \begin{equation}
 \begin{split}
     \mathcal{J}^{c}_{i}(\pi;x)&:=\mathbb{E}\bigg[\int_{t_i}^T L(s,X^\pi_{s},\pi_s)ds + g(X_T^{\pi})\Big|X^{\pi}_{t_i}=x\bigg],\\
     \hat{\mathcal{J}}_{i}(\pi;x)&:=\mathbb{E}\bigg[\sum_{j=i}^{N-1} L(t_{j},\hat{X}^\pi_{t_{j}},\pi_{t_{j}})\delta + g(\hat{X}_T^{\pi})\Big|\hat{X}^{\pi}_{t_i}=x\bigg].
 \end{split}
 \end{equation}
 We bound $|\hat{\mathcal{J}}_{i}(\pi;x)-\mathcal{J}^{c}_{i}(\pi;x)|$ in two steps, first replacing the continuous-time integral in $\mathcal{J}^c_i$ by its Riemann sum along the \emph{true} trajectory $X^\pi$, and then replacing $X^\pi$ by $\hat{X}^\pi$ at the grid points.

 By Assumption \textbf{A.2}, for any $\pi\in\Uppi_\delta$ and $j\in[N]$,
 \begin{equation}
 \begin{split}
     \bigg|\int_{t_{j}}^{t_{j+1}} \!\!\!\!\!\!\! L(s,{X}^\pi_{s},\pi_{t_j})ds-L(t_{j},{X}^\pi_{t_{j}},\pi_{t_j})\delta\bigg|
     &\le \int_{t_{j}}^{t_{j+1}} \!\!\!\!\!\!\!\big|L(s,{X}^\pi_{s},\pi_{t_j})-L(t_{j},{X}^\pi_{t_{j}},\pi_{t_j})\big|ds\\
     &\le C \int_{t_{j}}^{t_{j+1}}\!\!\!\!\!\!\!\!\!\!\!\! \big(\sqrt{s-t_j}+|{X}^\pi_{s}-{X}^\pi_{t_{j}}|\big)ds,
 \end{split}
 \end{equation}
 and, using $\mathbb{E}[|{X}^\pi_{s}-{X}^\pi_{t_{j}}|]\le C\sqrt{s-t_j}$, summing over $j\in\{i,\ldots,N-1\}$ and integrating gives
 \begin{equation}\label{eqn:J-Jgrid}
     \bigg|\mathcal{J}^{c}_{i}(\pi;x)-\mathbb{E}\bigg[\sum_{j=i}^{N-1} L(t_{j},{X}^\pi_{t_{j}},\pi_{t_{j}})\delta+g(X_T^\pi)\bigg]\bigg|\le C\sqrt{\delta}.
 \end{equation}
 This step never compares $\pi$ at two different arguments, only that $\pi$ is constant on each $[t_j,t_{j+1})$, so no regularity of $\pi$ itself is used.

 By \cite[Theorem~9.6.2]{KP13} and the Lipschitz continuity of $\mu,\sigma$ in $x$ (Assumption \textbf{A.2}), $\sup_{j\in[N+1]}\mathbb{E}[|{X}^\pi_{t_j}-\hat{X}^\pi_{t_j}|]\le C\sqrt{\delta}$, uniformly over $\pi\in\Uppi_\delta$. Together with the Lipschitz continuity of $L$ and $g$ in $x$ (Assumption \textbf{A.2}) and $\sum_{j=i}^{N-1} \delta\le T$, this gives
 \begin{equation}\label{eqn:Jgrid-hatJ}
     \bigg|\mathbb{E}\bigg[\sum_{j=i}^{N-1} L(t_{j},{X}^\pi_{t_{j}},\pi_{t_{j}})\delta+g(X_T^\pi)\bigg]-\hat{\mathcal{J}}_{i}(\pi;x)\bigg|\le C\sum_{j=i}^{N-1} \delta\sqrt{\delta}+C\sqrt{\delta}\le C\sqrt{\delta}.
 \end{equation}
 Combining \eqref{eqn:J-Jgrid} and \eqref{eqn:Jgrid-hatJ} by the triangle inequality,
 \begin{equation}\label{eqn:hatJ-J}
     |\hat{\mathcal{J}}_{i}(\pi;x)-\mathcal{J}^{c}_{i}(\pi;x)|\le C\sqrt{\delta},
 \end{equation}
 with $C$ depending only on the constants in Assumption \textbf{A.2}, uniformly over $\pi\in\Uppi_\delta$.

 Define $U^*(t_i,x):=\inf_{\pi\in\Uppi_\delta}\mathcal{J}^{c}_{i}(\pi;x)$, the value of the continuous-time problem restricted to controls piecewise-constant on the time grid. Since $\Uppi_\delta\subseteq\Uppi$, $U^*(t_i,x)\ge V(t_i,x)$. By Assumptions \textbf{A.1}--\textbf{A.2} (see the remark above), \cite[Theorem~2.1]{JPR19}, applied to the cost-minimization form of the problem, gives
 \begin{equation}\label{eqn:U*-V}
     0\le U^*(t_{i},x)-V(t_{i},x)\le C\delta^{1/4},
 \end{equation}
 with $C$ depending only on the constants in Assumptions \textbf{A.1}--\textbf{A.2}.

 Fix $\epsilon>0$ and choose $\pi^\epsilon\in\Uppi_\delta$ with $\mathcal{J}^{c}_{i}(\pi^\epsilon;x)\le U^*(t_i,x)+\epsilon$. By \eqref{eqn:hatJ-J} and \eqref{eqn:U*-V},
 \begin{equation}
     \hat{\mathcal{J}}_{i}(\pi^\epsilon;x)\le \mathcal{J}^{c}_{i}(\pi^\epsilon;x)+C\sqrt{\delta}\le V(t_i,x)+C\delta^{1/4}+C\sqrt{\delta}+\epsilon.
 \end{equation}
 Since $\pi^\epsilon\in\Uppi_\delta$, its restriction to the time grid is, by the identification of $\Uppi_\delta$ with the admissible class of \eqref{prob:control_discretized}, admissible for \eqref{prob:control_discretized}; since $\hat{V}$ is the infimum of $\hat{\mathcal{J}}_i(\cdot;x)$ over such admissible controls (see \eqref{defn:value_discrete}), $\hat{\mathcal{J}}_{i}(\pi^\epsilon;x)\ge \hat{V}(t_i,x)$. Letting $\epsilon\downarrow0$ and using $\sqrt{\delta}\le\delta^{1/4}$ for $\delta\le1$, we obtain
 \begin{equation}\label{eqn:half1}
     \sup_{x,i\in[N]}\big(\hat{V}(t_{i},x)-V(t_i,x)\big)\le C\delta^{1/4}.
 \end{equation}
To show $\sup_{x,i\in[N]}\big(V(t_{i},x)-\hat{V}(t_{i},x)\big)\le C\sqrt{\delta}$, fix $\epsilon>0$. Since $\hat{V}(t_i,x)$ is, by definition, the infimum of the discrete-time cost over admissible controls for \eqref{prob:control_discretized}, there exists an admissible control, $\pi^{n,\epsilon}$, for \eqref{prob:control_discretized} with
\begin{equation}\label{eqn:eps-opt-discrete}
    \mathbb{E}\bigg[\sum_{j=i}^{N-1} L(t_{j},\hat{X}^{n,\epsilon}_{t_{j}},{\pi}^{n,\epsilon}_{t_{j}})\delta + g(\hat{X}_T^{n,\epsilon})\Big|\hat{X}^{n,\epsilon}_{t_i}=x\bigg]\le \hat{V}(t_i,x)+\epsilon,
\end{equation}
where $\hat{X}^{n,\epsilon}$ is the corresponding trajectory of \eqref{prob:control_discretized}. We construct from $\pi^{n,\epsilon}$ an admissible continuous-time control $\pi$, and compare its continuous-time cost $U(t_i,x)$ on one side to $V(t_i,x)$ -- by admissibility -- and on the other to the left-hand side of \eqref{eqn:eps-opt-discrete} -- by the same quadrature-and-discretization argument used above -- giving $V(t_i,x)\le U(t_i,x)\le \hat{V}(t_i,x)+\epsilon+C\sqrt{\delta}$; letting $\epsilon\downarrow0$ then gives the claim.

Set $\pi_t:=\pi^{n,\epsilon}_{t_j}$ for $t\in[t_j,t_{j+1})$ and $j\in\{i,\ldots,N-1\}$. By construction, $\hat{X}^{n,\epsilon}$ is the Euler-Maruyama discretization of the continuous-time process $X^\pi$ solving
\begin{equation}
    dX_t^\pi =\mu(t,X_t^{\pi},\pi_t)dt+\sigma(t,X_t^{\pi},\pi_t)dW_t,
\end{equation}
so by \cite[Theorem~9.6.2]{KP13} and Assumption~\textbf{A.2}, $\mathbb{E}[|X^{\pi}_{t_j}-\hat{X}^{n,\epsilon}_{t_j}|]\le C\sqrt{\delta}$; as before, $\mathbb{E}[|X^\pi_{s}-X^\pi_{t_j}|]\le C\sqrt{s-t_j}$. Define
\begin{equation}
    {U}(t_i,x):=\mathbb{E}\bigg[\int_{t_i}^T L(s,{X}^{\pi}_{s},{\pi}_s)ds + g({X}^{\pi}_T)\bigg|{X}^\pi_{t_i}=x\bigg].
\end{equation}
Since $\pi$ is an admissible control for the continuous-time problem \eqref{prob:control} started from $x$ at $t_i$, and $V(t_i,x)$ is the infimum of the cost over all such controls, $U(t_i,x)\ge V(t_i,x)$.

It remains to bound $U(t_i,x)$ by the left-hand side of \eqref{eqn:eps-opt-discrete}. Since $\pi_s=\pi_{t_j}$ for $s\in[t_j,t_{j+1})$, the triangle inequality gives
\begin{equation*}
\begin{split}
    &\bigg|\int_{t_{j}}^{t_{j+1}} \!\!\!\!\!\!\! L(s,{X}^{\pi}_{s},{\pi}_{t_j})ds-L(t_{j},\hat{X}^{n,\epsilon}_{t_{j}},{\pi}_{t_j})\delta\bigg|
    \le \!\!\! \int_{t_{j}}^{t_{j+1}} \!\!\!\!\!\!\!\big|L(s,{X}^\pi_{s},{\pi}_{t_j})-L(t_{j},\hat{X}^{n,\epsilon}_{t_{j}},{\pi}_{t_j})\big|ds\\
    &\le \!\!\! \int_{t_{j}}^{t_{j+1}} \!\!\!\!\!\!\!\big|L(s,{X}^\pi_{s},{\pi}_{t_j})-L(s,\hat{X}^{n,\epsilon}_{t_{j}},{\pi}_{t_j})\big|ds+\!\!\!\int_{t_{j}}^{t_{j+1}} \!\!\!\!\!\!\!\big|L(s,\hat{X}^{n,\epsilon}_{t_{j}},{\pi}_{t_j})-L(t_{j},\hat{X}^{n,\epsilon}_{t_{j}},{\pi}_{t_j})\big|ds
\end{split}
\end{equation*}
By the $\sfrac{1}{2}$-Hölder continuity of $L$ in $t$ (Assumption~\textbf{A.2}),
\begin{equation}\label{Ls-Lt_j}
    \big|L(s,\hat{X}^{n,\epsilon}_{t_{j}},{\pi}_{t_j})-L(t_{j},\hat{X}^{n,\epsilon}_{t_{j}},{\pi}_{t_j})\big|\le C \sqrt{s-t_j},
\end{equation}
and, inserting $X^\pi_{t_j}$ as an intermediate point, the Lipschitz continuity of $L$ in $x$ (Assumption~\textbf{A.2}) gives
\begin{equation}\label{X-hatXn}
    \big|L(s,{X}^\pi_{s},{\pi}_{t_j})-L(s,\hat{X}^{n,\epsilon}_{t_{j}},{\pi}_{t_j})\big|
    \le C \big(|{X}^\pi_{s}-{X}^{\pi}_{t_{j}}|+ |{X}^\pi_{t_j}-\hat{X}^{n,\epsilon}_{t_{j}}|\big)
\end{equation}
Combining \eqref{Ls-Lt_j}--\eqref{X-hatXn}, summing over $j\in\{i,\ldots,N-1\}$, integrating, and using $\mathbb{E}[|X^{\pi}_{t_j}-\hat{X}^{n,\epsilon}_{t_j}|]\le C\sqrt{\delta}$ and $\mathbb{E}[|{X}^\pi_{s}-{X}^{\pi}_{t_{j}}|]\le C\sqrt{s-t_j}$, we obtain
\begin{equation*}
\begin{split}
    U(t_i,x)&\le \mathbb{E}\bigg[\sum_{j=i}^{N-1} L(t_{j},\hat{X}^{n,\epsilon}_{t_{j}},{\pi}^{n,\epsilon}_{t_{j}})\delta + g(\hat{X}_T^{n,\epsilon})\Big|\hat{X}^{n,\epsilon}_{t_i}=x\bigg]+C\sqrt{\delta}\\
    &\le \hat{V}(t_i,x)+\epsilon+C\sqrt{\delta},
\end{split}
\end{equation*}
the last step by \eqref{eqn:eps-opt-discrete}. Since $V(t_i,x)\le U(t_i,x)$ and $\epsilon>0$ was arbitrary,
\begin{equation*}
    \sup_{x,i\in[N]}\big(V(t_{i},x)-\hat{V}(t_{i},x)\big)\le C\sqrt{\delta}.
\end{equation*}
 \end{proof}
 \begin{proof}[Proof of Proposition~\ref{prop:epsilon_optimal}]
It is sufficient to show that if $\sup_x|V(t_m,x)-\hat{V}(t_m,x)|\le \dfrac{\epsilon}{2}$, then
\begin{equation}\label{A20}
     0 \le\mathbb{E}\bigg[\int_{t_n}^{t_m}L(r,X^{\pi^\epsilon}_r,\pi^\epsilon_r)dr + V(t_m,X_{t_m}^\pi)\Big|X_{t_n}=x\bigg]-V(t_n,x)\le {\epsilon}
\end{equation}
The first inequality above is the direct result of \eqref{DPP_continuous}. Fix $\eta>0$; since $V(t_n,x)$ is, by \eqref{DPP_continuous}, the infimum of the right-hand side over $\pi\in\Uppi_{t_n,t_m}$, no assumption is needed to guarantee the existence of an admissible $\pi^\eta\in\Uppi_{t_n,t_m}$ with
\begin{equation}
    \mathbb{E}\bigg[\int_{t_n}^{t_m}L(r,X_r^{\pi^\eta},\pi^\eta_r)dr+V(t_m,X_{t_m}^{\pi^\eta})\Big|X_{t_n}=x\bigg]\le V(t_n,x)+\eta.
\end{equation}
Let $\hat{X}^{\pi^\eta}$ be the Euler-Maruyama discretization of ${X}^{\pi^\eta}$. Therefore,
 for $r\in[t_j,t_{j+1})$, $\mathbb{E}[|{X}^{\pi^\eta}_{r}-\hat{X}^{\pi^\eta}_{r}|]\le C\sqrt{\delta}$ and $\mathbb{E}[|{X}^{\pi^\eta}_{r}-{X}^{\pi^\eta}_{t_j}|]\le C\sqrt{r-t_j}$. In particular,
\begin{equation}
    \begin{split}
        &\bigg|\mathbb{E} \bigg[\int_{t_{j}}^{t_{j+1}}\big(L(r,X^{\pi^\eta}_r,\pi^\eta_{r})-L(t_{j},\hat{X}^{\pi^\eta}_{t_j},{\pi}^\eta_{t_j})\big)dr\bigg]\bigg|\\
        &\le C\mathbb{E} \bigg[\int_{t_{j}}^{t_{j+1}}\big(\sqrt{r-t_j}+|X^{\pi^\eta}_{r}-X^{\pi^\eta}_{t_j}|\big)dr\bigg]\le C\delta^{3/2}
    \end{split}
\end{equation}
Since the Euler-Maruyama discretization of ${X}^{\pi^\epsilon}$ is the same as $\hat{X}^{*n}$, we have
\begin{equation}
    \begin{split}
        &\bigg|\mathbb{E} \bigg[\int_{t_{j}}^{t_{j+1}}\big(L(r,X^{\pi^\epsilon}_r,\pi^\epsilon_r)-L(t_j,\hat{X}^{*n}_{t_j},{\pi}^*_{t_j})\big)dr\bigg]\bigg|\\
        &\le C\mathbb{E} \bigg[\int_{t_{j}}^{t_{j+1}}\big(\sqrt{r-t_j}+|X^{\pi^{\epsilon}}_{r}-X^{\pi^{\epsilon}}_{t_j}|\big)dr\bigg]\le C\delta^{3/2}
    \end{split}
\end{equation}
One can use the above estimates to write:
\begin{equation}
    \begin{split}
        V(t_n,x)+\eta &\ge \mathbb{E}\bigg[\int_{t_n}^{t_m}L(r,X^{\pi^\eta}_r,\pi^\eta_r)dr + V(t_m,X_{t_m}^{\pi^\eta})\bigg]\\
        &\ge - C\sqrt{\delta} -\frac{\epsilon}{2} + \mathbb{E}\bigg[\int_{t_n}^{t_m}L(t_j,\hat{X}^{\pi^\eta}_{t_j},{\pi}^{\eta}_{t_j})dr+\hat{V}(t_m,\hat{X}_{t_m}^{\pi^\eta})\bigg]\\
        &\ge - C\sqrt{\delta} -\frac{\epsilon}{2} + \mathbb{E}\bigg[\int_{t_n}^{t_m}L(t_j,\hat{X}^{\pi^\epsilon}_{t_j},{\pi}^{\epsilon}_{t_j})dr+\hat{V}(t_m,\hat{X}_{t_m}^{\pi^\epsilon})\bigg]
    \end{split}
\end{equation}
Letting $\eta\downarrow0$ and using $C\sqrt\delta\le \epsilon/2$, \eqref{A20} holds true.
\end{proof}

\bibliographystyle{siamplain}
\bibliography{refs}
\end{document}